\documentclass[12pt]{article}

\usepackage{amsmath,amssymb,mathrsfs,theorem,subfiles,bm}
\usepackage{graphicx}
\usepackage{subfigure}
\usepackage[top=30truemm, bottom=30truemm,left=25truemm,right=25truemm]{geometry}
\usepackage{tikz}  
\usetikzlibrary{shapes.geometric, arrows} 
\usepackage{capt-of}
\usepackage[thinlines]{easytable}
\usepackage{pbox}
\usepackage{comment}
\usepackage{setspace}
\usepackage{marvosym}
\usepackage[colorlinks,allcolors=blue]{hyperref}  

\newcommand{\tx}{x}
\newcommand{\tP}{P}

\newcommand{\ty}{y}

\def\e{{\mathrm{e}}}

\def\t0{0}

\def\tX{\textbf{X}}
\def\tK{\textbf{K}}
\def\tR{\textbf{R}}
\def\ia{I_{a}}
\def\is{I_{s}}
\def\ra{R_{a}}
\def\rs{R_{s}}
\def\rh{R}
\def\dh{D}

\def\iia{i_{a}}
\def\iis{i_{s}}
\def\rra{r_{a}}
\def\rrs{r_{s}}
\def\rrh{r}
\def\ddh{d}
\def\ba{\beta_{a}}
\def\bs{\beta_{s}}
\def\ga{\gamma_{a}}
\def\gs{\gamma_{s}}
\def\th{\tau_{H}}
\def\gh{\gamma_{H}}
\def\gp{\gamma_{D}}  
\def\rt{R_{t}}
\def\ro{R_{0}}
\def\N{\mathbb N}

\title{Analysis of COVID-19 in Japan with Extended SEIR model and ensemble Kalman filter}

\begin{document}
\tikzstyle{block} = [rectangle, draw, thick = 1,
    text width=2.5em, text centered, rounded corners, minimum height=2.5em]
\tikzstyle{block2} = [rectangle, draw, thick = 1,
    text width=5.5em, text centered, rounded corners, minimum height=2.5em]
\tikzstyle{block1} = [trapezium, trapezium left angle = 65,
    trapezium right angle = 115, draw, thick = 1,
    text width=3em, text centered, rounded corners, minimum height=2.5em]
\tikzstyle{line} = [draw, thick = 1, -latex']
\tikzstyle{cloud} = [draw, ellipse,fill=red!20, node distance=3cm,
    minimum height=2em]

\author{Q. Sun${}^{a,b}$\Letter, T. Miyoshi${}^{a,c,d}$, S. Richard${}^{a,b}$}

\date{\small}
\maketitle
\vspace{-1cm}

\begin{quote}
\emph{
\begin{enumerate}
\item[a)] Data Assimilation Research Team, RIKEN Center for Computational Science (R-CCS), Kobe, Japan
\item[b)] Graduate School of Mathematics, Nagoya University, Nagoya, Japan
\item[c)] Prediction Science Laboratory, RIKEN Cluster for Pioneering Research, Kobe, Japan
\item[d)] RIKEN interdisciplinary Theoretical and Mathematical Sciences Program (iTHEMS), Kobe, Japan
\item[] \emph{E-mails:} qiwen.sun@riken.jp, takemasa.miyoshi@riken.jp, \\
richard@math.nagoya-u.ac.jp
\end{enumerate}
}
\end{quote}

\begin{abstract}
{
We introduce an extended SEIR infectious disease model with data assimilation 
for the study of the spread of COVID-19.
In this framework, undetected asymptomatic and pre-symptomatic cases are taken into account,
and the impact of their uncertain proportion is fully investigated. 
The standard SEIR model does not consider these populations, while their role in the propagation 
of the disease is acknowledged.
An ensemble Kalman filter is implemented to assimilate reliable observations of three compartments in the model. The system tracks the evolution of the effective reproduction number and estimates the unobservable subpopulations. The analysis is carried out for three main prefectures of Japan and for the entire population of Japan. 
For these four populations, our estimated effective reproduction numbers are more stable than
the corresponding ones estimated by a different method (Toyokeizai).
We also perform sensitivity tests for different values of some uncertain medical parameters, like the relative infectivity of symptomatic / asymptomatic cases. 
The regional analysis results suggest the decreasing efficiency of the states of emergency.
}
\end{abstract}

\textbf{keywords:} COVID-19, data assimilation, ensemble Kalman filter, extended SEIR model


\section{Introduction}\label{sec:int}

The outbreak of coronavirus 2019 (COVID-19) had a huge impact on human society,
and its devastating effects are still present more than 18 months after its official acknowledgment.
One of the specific characteristics of this disease increases the difficulties for any scientific research based on infectious disease models, namely the existence
of asymptomatic cases.
Their existence, but also their potential infectiousness, blur the full picture of the disease spreading inside the society \cite{WHO}. 
The simultaneous existence of mild symptomatic cases, which often remain undocumented,
also generates an additional challenge for any epidemiological investigations.
For example, one study discussing the potential fraction of undocumented cases at the beginning of the outbreak ends up with a ratio of $86\%$  undocumented infections \cite{LI}. 

There exists now a huge literature related to COVID-19, but only a small number of papers are using techniques of data assimilation, as we shall do. For that reason, we present only a few references closely related to our investigations. 
To get a better understanding of the spread of COVID-19, 
some studies start with a specially designed infectious disease model with appropriate structures. Those elaborated models may accurately describe the process in reality but they also bring more uncertainties and unknown parameters.
For example in \cite{Arm}, a $22$ variables model is  adopted with each variable representing a subpopulation. 
The structure related to asymptomatic, pre-symptomatic, detected, and undetected cases are all properly described by the structure. 
However, the inference is done with only observations of $4$ compartments among all variables.
Another example is \cite{E}, where an extended SEIR model with $11$ age-classes is used to estimate the posterior distribution of parameters and make short period predictions. The distribution of mild, severe and fatal cases for each age class is introduced as pre-defined parameters. 
Also, a matrix which describes the relative transmissions between each pair of age groups is  predefined. 

Despite these elaborated models, most of the investigations are still performed on the common infectious disease models, namely the SIR and SEIR models.
Such models are easier to implement and contain less parameters, but they also sacrifice some details of the process.
For example the SEIR model is used to study the dynamical behavior of COVID-19 outbreaks at the regional level in \cite{EN},
and the SIR model was used to track the effective reproduction number for $124$ countries in \cite{AM}.
A slightly more elaborated model is also introduced in \cite{NA}, namely the SITR model with T standing for \emph{under treatment}, for inferring infection numbers and parameters values. 

For dealing with real data, one of the recent tools employed for the study of 
infectious diseases is data assimilation. 
The techniques developed for data assimilation have the ability to optimally meld dynamical systems with noisy and imperfect observations. They also provide forecasts and estimations of parameters and variables \cite{RH}. 
As examples of investigations studying the details of the implementation of these techniques to epidemic models, let us mention 
\cite{MI} and \cite{RE}. 
The previously cited work \cite{Arm} also uses statistical data assimilation for identifying the measurements required for getting accurate states and parameters estimations. 
Finally, in \cite{GH} an ensemble Kalman filter is adopted to forecast the COVID-19 pandemic in Saudi Arabia with an extended SEIR model including vaccination.  

Let us now present our research. First, we construct an extended SEIR model which takes into account one of the main specificities of COVID-19: the simultaneous existence of an asymptomatic\;\!/\;\!pre-symptomatic population and of a symptomatic population. These two infected populations have different characteristics which are encoded with different parameters.
The model is also constructed such that it can be fed with the data of only three compartments, namely the Hospitalized (or treated) population, the Recovered population, and the Deceased population, naturally identified by the letters $H$, $R$, and $D$, respectively.
Note that the population $R$ or $D$ are coming from the population $H$, and these three populations are considered as fully recorded by health services. Except for the first couple of weeks at the beginning of the epidemic, these data are also considered as the most reliable ones.
With this model and these data, our main aim is to provide information about 
the effective reproduction number and about the population of asymptomatic\;\!/\;\!pre-symptomatic
or undetected symptomatic, with uncertainty ranges. 
For these investigations, we shall also mainly concentrate on the population of Tokyo, 
but provide additional analysis with two additional prefectures (Osaka and Kanagawa) and with the entire population of Japan.
Let us mention that targeting a specific population gives the opportunity to rely on local parameters and to use information provided by local medical research or health services.
We also emphasize that for data assimilation, we use
the ensemble transform Kalman filter \cite{HKS} which is efficient compared with 
the standard extended Kalman filter.
 
As a result of our investigations we get an effective reproduction number for the four populations mentioned before, as well as a time dependent estimate for the number of asymptomatic\;\!/\;\!pre-symptomatic in these four populations.  
A comparison between our estimated effective reproduction number and the same statistic provided by Toyokeizai \cite{TK} shows the reliability of our strategy.
As an asset of our approach, we test the sensitivity of our model and its outcomes to the values of  uncertain parameters borrowed from the literature.  
These experiment's results show a robust performance of the strategy used in our investigations. Another clear result from our computation of the effective reproduction number is the decay of the effectiveness of the states of emergency as their number increases. 
This effect is clearly visible independently for the three prefectures and for the entire country. 
Note that a similar effect is also visible in the use of public transportation, see for example
\cite{MLIT}.
Other outcomes of our investigation are gathered in Sections  \ref{sec:exp} and \ref{sec:disc}.

Let us finally describe the structure of this paper. In Section \ref{sec:esm}, 
we recall the standard SIR and SEIR models, 
and introduce the extended SEIR model. We also provide information for the computation of the effective reproduction number, and discuss the constant parameter values and the observations. 
In Section \ref{enkf}, we explain the technical details of ensemble transform Kalman filter. Readers familiar with data assimilation and Kalman filters can skip this section without any loss for the application to infectious diseases. 
In Section \ref{sec:exp}, background information about the experiments, and subsequently technical information, are provided. The main results of our investigations are also presented in this section. 
The discussion and the comparisons between the different populations are provided in
Section \ref{sec:disc}. We finally prepare the ground for future projects.

\section{Extended SEIR model}\label{sec:esm}

In this section, we recall a few information on the SIR and SEIR models, and introduce the extended SEIR model.
We also discuss medical parameters, and provide some explanations about the observations.

\subsection{The SIR  and SEIR models}

The SIR model is a deterministic epidemic model expressed by a system of differential equations. 
Its construction is based on Kermack–McKendrick theory, and it describes the transmission process of an infectious disease. 
A given population of size $N$ is divided into three mutually exclusive sub-populations call compartments: the susceptible hosts $S$, the infectious hosts $I$, and the recovered hosts $R$. 
For any given time $t$, $S(t)$, $I(t)$, and $R(t)$ describe the number of individuals in each respective compartment, and they satisfy $S(t)+I(t)+R(t)=N$.
Individuals in compartment $S$ can be infected by individuals in compartment $I$ through direct contact. 
Once a successful contact happens, the infected individual leaves compartment $S$ and becomes a member of compartment $I$. 
The number of newly infected individuals per unit time is given by $\beta I(t)S(t)/N$, where the parameter $\beta$ is called the transmission coefficient. 
Similarly, once a patient recovers, this person leaves compartment $I$ and goes to compartment $R$ immediately. 
The transfer between $I$ and $R$ is controlled by the recovery rate $\gamma$.
The transfer rate, namely the number of transfer individuals per unit time, is given 
by $\gamma I(t)$.

\begin{center}
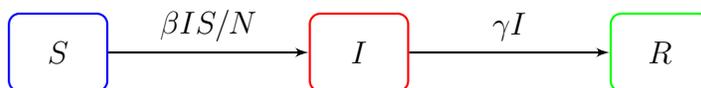

\begin{tikzpicture}[node distance=4cm, auto]
\node [block, draw = blue] (init) {$S$};
\node [block, draw = red, right of=init] (I) {$I$};
\node [block, draw = green, right of= I] (R) {$R$};
\path [line] (init) -- node {$\beta IS/N$}(I);
\path [line] (I) -- node {$\gamma I$}(R);
\end{tikzpicture}
\captionof{figure}{Transfer diagram for the SIR model}
\label{psir}
\end{center}

Figure \ref{psir} shows the process of the SIR model, and the following differential system describes its dynamic:
\begin{align}
\frac{dS}{dt} &= -\beta IS/N \nonumber,\\
\frac{dI}{dt} &= \beta IS/N - \gamma I, \label{sirm}\\
\frac{dR}{dt} &=\gamma I. \nonumber
\end{align}
Note that the system \eqref{sirm} assumes a permanent immunity once recovered.
In other words, no individual can be infected a second time.  
The model also assumes a constant total population: there is no inflow to the system, or outflow from the system. 

The SEIR model is very similar to the previous model, but with an additional compartment $E$ between the compartments $S$ and $I$. The individuals in $E$ are exposed, namely they have been contaminated, but they are not infectious yet. The time spent in $E$ corresponds to the incubation period, before becoming infectious.
The transfer rate from $E$ to $I$ is given by $\varepsilon E(t)$.

\begin{center}
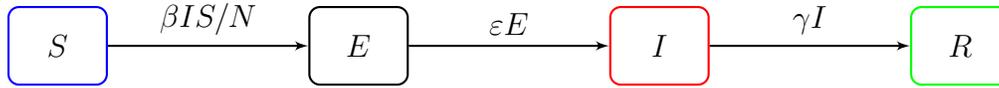

\begin{tikzpicture}[node distance=4cm, auto]
\node [block, draw = blue] (init) {$S$};
\node [block, draw = black, right of=init] (E) {$E$};
\node [block, draw = red, right of=E] (I) {$I$};
\node [block, draw = green, right of= I] (R) {$R$};
\path [line] (init) -- node {$\beta IS/N$}(E);
\path [line] (E) -- node {$\varepsilon E$}(I);
\path [line] (I) -- node {$\gamma I$}(R);
\end{tikzpicture}
\captionof{figure}{Transfer diagram for the SEIR model}
\label{pseir}
\end{center}

Figure \ref{pseir} shows the process of the SEIR model, and the following differential system describes its dynamic:
\begin{align*}
\frac{dS}{dt} &= -\beta IS/N,\\
\frac{dE}{dt} &= \beta IS/N - \varepsilon E, \\
\frac{dI}{dt} &= \varepsilon E- \gamma I, \\
\frac{dR}{dt} &=\gamma I. 
\end{align*}

\subsection{The extended SEIR model}

The extended SEIR model has been developed based on some specific features
of the COVID-19 epidemic, as described now. 
In early 2020, health services already noticed that some infected individuals, who did
not show any symptom, could spread the disease. 
These persons correspond either to asymptomatic hosts or to pre-symptomatic
hosts.
In the first cohort, people will never show any symptom, while the second
cohort corresponds to individuals who will show some symptoms subsequently.  
Clearly, asymptomatic hosts and pre-symptomatic hosts are difficult to be identified by health services, even though they play an important role in the spread of the disease.
Symptomatic individuals are also infectious, but they can be more easily identified 
precisely because of their symptoms. Thus, if symptomatic individuals correspond to the compartment $I$ of the above SIR or SEIR models, then there is no compartment left for the asymptomatic or the pre-symptomatic individuals.  

The existence of COVID-19 infectious spreaders who do not show (yet or at all) any symptom already brings a lot of uncertainties to health services. 
In addition, whenever the symptoms are relatively mild, 
some symptomatic individuals do hesitate to report the health service 
\cite[p.~11]{OP}. As a consequence, daily new cases are under-reported.
Combining this effect with some delayed information, some inaccurate tests (especially
in the first phase of the epidemic), and other reasons that we are not aware of,
it turns out that that the reported data are not very accurate. In such a situation one needs to use proper strategies to analyse the observation data. 

\begin{center}
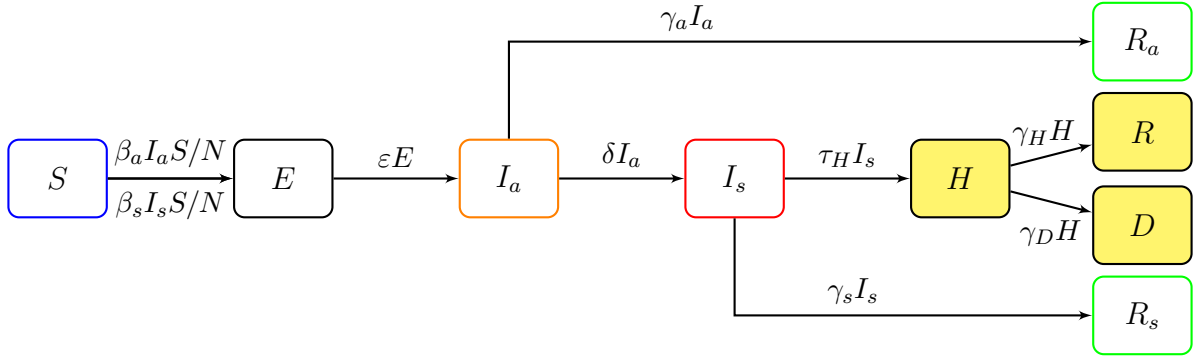

\begin{tikzpicture}[node distance=3cm, auto]
\node [block, draw = blue] (init) {$S$};
\node [block, draw = black, right of=init] (E) {$E$};
\node [block, draw = orange, right of=E] (I1) {$\ia$};
\node [block, draw = red, right of=I1] (I2) {$\is$};
\node [block, fill=yellow!70, right of= I2] (H) {$H$};
\node [block,  fill=yellow!70, above right of= H, yshift=-1.5cm, xshift = 0.3cm] (RH) {$\rh$};
\node [block,  fill=yellow!70, below right of= H, yshift= 1.5cm, xshift = 0.3cm] (DH) {$\dh$};
\node [block, draw = green, above of = RH, yshift=-1.8cm] (R1) {$\ra$};
\node [block, draw = green, below of = DH, yshift=1.8cm] (R2) {$\rs$};
\path [line] (init) -- node {\small $\ba \ia S/N$}(E);
\path [line] (init) -- node [shift={(0,-7mm)}] {\small  $\bs \is S/N$}(E);
\path [line] (E) -- node {\small $\varepsilon E$}(I1);
\path [line] (I1) -- node {\small $\delta \ia $}(I2);
\path [line] (I2) -- node {\small $\tau_{H} \is $}(H);
\path [line] (H) -- node [near end, shift={(2mm,-1mm)}] {\small $\gh H$}(RH);
\path [line] (H) -- node [near start, shift={(-3mm,-8mm)}] {\small $\gp H$}(DH);
\path [line] (I1) |- node [near end, xshift = -1.5cm] {\small $\ga \ia$}(R1);
\path [line] (I2) |- node [near end, xshift = -0.8cm] {\small $\gs \is $}(R2);
\end{tikzpicture}
\captionof{figure}{Transfer diagram for the extended SEIR model. Compartment $I$ of the SEIR model is divided into two compartments $\ia$ (asymptomatic\;\!/\;\!pre-symptomatic) and $\is$ (symptomatic)} 
\label{pesir}
\end{center}

In order to meet the special features of COVID-19, we separate the compartment $I$ of the SEIR model into two compartments $\ia$ and $\is$. 
These compartments correspond to  asymptomatic\;\!/\;\!pre-symptomatic and to symptomatic states, respectively. Both $\ia$ and $\is$ can infect $S$. 
As shown in Figure \ref{pesir}, the transmission processes related to $\ia$ and to $\is$ have transmission coefficients $\ba$ and $\bs$ respectively. 
Once an individual in $S$ gets infected, this person becomes a member of $E$,
and then moves to $\ia$ when it becomes infectious.
In $\ia$, part of these individuals (the asymptomatic) will never generate any symptoms. They will thus spend some time in $\ia$, and then recover. The corresponding compartment is denoted by $\ra$. In contrast, the other part of individuals in $\ia$ will develop symptoms, and then move to $\is$. 
In $\is$, a majority of individuals will be identified by health services, but as mentioned above a minority will recover without being identified by any health services. 
Compartment $\rs$ corresponds to those recovered individuals who have not been registered. 
The identified persons in $\is$ will then move to the compartment $H$,
which corresponds to hospitalized or treated patients. The ones staying at home but under
medical control, or the ones isolated in a hotel, are all included in the compartment
$H$.
Finally, the patients in $H$ will recover, and thus move to the compartment $\rh$, or will decease and end up in the compartment $\dh$. 
As for the SIR or SEIR models, we assume a permanent immunity, which means that there is no flow from $\ra$, $\rs$, or $\rh$ to $S$. Also, the total number $N$ of individual is constant, namely at all time one has
\begin{equation}\label{eq_conservation}
N = S + E + \ia +\is + H + \rh + \dh + \ra + \rs.
\end{equation}

Based on the Figure \ref{pesir} and on the explanations provided above, the differential system of the extended SEIR model reads:
\begin{align}\label{esirm}
\frac{dS}{dt} &= -\ba\ia S/N - \bs\is S/N   &   \frac{d \rh}{dt} &=  \gh H \nonumber \\
\frac{dE}{dt} & =  \ba\ia S/N  + \bs\is S/N - \varepsilon E &  \frac{d \dh}{dt} &= \gp H \nonumber  \\
\frac{d \ia}{dt} &= \varepsilon E - \delta \ia - \ga \ia     &       \frac{d \ra}{dt} &=  \ga \ia  \\
\frac{d \is}{dt} &= \delta \ia - \tau_{H} \is - \gs \is &  \frac{d \rs}{dt} &=  \gs \is   \nonumber\\    
\frac{dH}{dt} &=  \tau_{H} \is - \gh H - \gp H      &  {} \nonumber 
\end{align}
where $\ba ,\bs, \varepsilon, \delta, \ga, \gs, \gp$ and $\gh$
are some medical parameters. Note that some of them are time dependent.

\subsection{The reproduction number}

The basic reproduction number, denoted by $\ro$, can be interpreted as the average number of secondary cases generated by on primary case in a susceptible population.
It is commonly admitted that $\ro = 1$ is a threshold value, as explained below. 
We also refer to \cite{DW,H} for more explanations and more precise statements. 

To study $\ro$ for different models, a general definition of the basic reproduction number is introduced based on the notion of disease free equilibrium (DFE).
In a DFE the population remains in the absence of disease. For example, in the SIR or SEIR models, it means that $I = 0$ while in the extended SEIR model it means that $\ia = \is = 0$.
In this context, the basic reproduction number can be defined as the number of new
infections produced by a typical infectious individual in a population at a DFE.
The main feature of $\ro$ is that it corresponds to a threshold parameter, namely if $\ro < 1$, the DFE is locally asymptotically stable, while if $\ro > 1$, the DFE is unstable and an outbreak is possible.

The precise expression for $\ro$ is clearly model dependent, but numerous examples
are available in the literature. 
For example, let us consider a general compartmental disease transmission model.
Such models are represented by a system of ordinary differential equations, and  
under natural assumptions it can be shown that these models have a DFE, see \cite{DW}.
In this reference, the precise expression for $\ro$ is then provided for some classes of models, and the extended SEIR model meets the assumptions of the staged progression model, as presented in \cite[Sec.~4.3]{DW}. 
For the extended SEIR model, the expression for $\ro$ then reads: 
\begin{align}
\ro =\frac{\ba}{\delta + \ga} + \frac{\bs \delta}{(\delta + \ga)(\gs + \tau_{H})}, \label{eqro}
\end{align}
where $\ba$ and $\bs$ are the initial transmission coefficients at time $0$.

To understand the above expression, observe that $\delta / (\delta + \ga)$ 
corresponds to the fraction of individuals of $\ia$ going to the compartment $\is$, while
$1/(\delta + \ga)$ and $1/(\gs + \th)$ define the average times an infected individual
spends in compartments $\ia$ and $\is$ respectively.
Thus, each term on the R.H.S.~of \eqref{eqro} represents the number of new infections generated by an infectious individual during the time spent in the compartments $\ia$ and $\is$.

In contrast, the definition of the effective reproduction number $\rt$ at time $t$ is much more delicate. The various challenges and possible errors have recently been  discussed in \cite{GMc}, to which we refer for a thorough discussion. 
In our setting, we shall keep the interpretation of $\rt$ as
the average number of secondary cases generated by one primary case.
This approach is possible because the transmission coefficients at time $t$ 
are available in our approach, and therefore one can compute $\rt$ with \eqref{eqro} and the corresponding $\ba$ and $\bs$ at time $t$. Additional information on $\rt$ will
be provided in Section \ref{sec:disc}.

\subsection{The medical parameters}\label{sec:mpara}

It clearly appears in Figure \ref{pesir} and in the corresponding system \eqref{esirm}
that several parameters control the flows between the compartments. 
The values of these parameters may result in very different behaviors of the model. 
We refer for example to \cite[Sec.~1.4]{L} and to \cite[Chap.~2]{BDW} for 
general discussions of model behaviors and the role of parameters.
In our setting, some parameters are easy to evaluate, as for example the recovery rate 
$\gh$ or the death rate $\gp$, but others are difficult to estimate, as for example the transmission coefficients $\ba$ and $\bs$.
Let us also stress that some parameters depend on location. 
Since our experiment is based on data from Japan, see Section \ref{sec:exp},
the medical parameters are chosen accordingly. 
For that reason, we use research or survey results specific to Japan, or even more 
precisely to specific prefectures in Japan. 

In Table \ref{tbpara} we list the estimated values of some parameters for the extended SEIR model, and provide the sources of information. 
Several parameters in the table have the form of the product of a percentage and the inverse of a time length. 
A similar setting for the construction of the parameters can be found for example in \cite{E,KI}. For $\delta$ and $\ga$, the percentage parts should sum up to $1$.
For the parameters $\th$ and $\gs$, some information can be found in the surveys \cite{OP, TBP} and the health services website \cite{TC}.
For parameter $\gh$ and $\gp$, instead of using constant value estimated by health services, we shall use the observation data to estimate their values on a daily basis. The details will be discussed in Section \ref{sec:exp}.

Let us stress that the value assigned to some of these parameters has evolved during 
the last 12 months. For example, the ratio of asymptomatic individuals was thought to be very high at the beginning of the epidemic (up to $80\%$), 
but some recent research \cite{BCB, BEC, PL} have revised this ratio to $17\%$ to $20\%$.
Our knowledge about the length of infectious periods has also evolved, and the possible values are spread over a rather long interval. In Table \ref{tbpara} we list some mean values, but in the simulations additional uncertainties are implemented.
Since pre-symptomatic patients become infected before the appearance of symptoms, the incubation period (encoded in $\varepsilon$) has been shorten a little bit, and the last 2 days of this incubation period have been moved to $\ia$. 

One very delicate question is the relation between $\ba$ and $\bs$, namely between the transmission coefficient for asymptomatic\;\!/\;\!pre-symptomatic and the transmission coefficient for symptomatic. For our investigations, we shall rely on the result of the systematic review \cite{BCB} which asserts that the relative risk of asymptomatic transmission
is $42\%$ lower than that for symptomatic transmission. As a consequence, we shall fix
\begin{equation}\label{eq_on_beta}
\ba= 0.58 \bs \quad \hbox{ or equivalently }\quad
\bs= 1.72 \ba. 
\end{equation} 
This factor $0.58$ is slightly smaller but of a comparable scale compared to earlier investigations, see for example \cite{BEC}.
\begin{table}[h]
\centering
\begin{TAB}(r,1cm,1.2cm)[2pt]{c|c|c|l}{|c|c|c|c|c|c|c|c|}
 parameter & estimation & source & \hfil remark \\
 $\varepsilon$ & $(3\textrm{ days})^{-1}$   & \cite{McA} &  incubation period, not contagious \\
 $\delta$ & \pbox{6cm}{$83\% \times( 2\textrm{ days})^{-1}$  \\
 $(95\% \hbox{ CI }  80\% \hbox{ to }86\%)$} & \cite{BCB, McA} & \pbox{15cm}{proportion of pre-symptomatic \\ $\times$ $(\textrm{duration of pre-symptomatic})^{-1}$} \\
 $\tau_{H}$ & \pbox{6cm}{$78\% \times (8.3 \textrm{ days})^{-1}$ ($\sim$ 2020/5/31)\\
$78\% \times (5.2 \textrm{ days})^{-1}$ (2020/6/1 $\sim$)}   & \cite{OP, M3} & \pbox{15cm}{proportion of detected symptomatic \\ $\times$ $(\textrm{days of symptoms onset})^{-1}$} \\
$\ga$ & \pbox{6cm}{$17\% \times (9 \textrm{ days})^{-1}$ \\
$(95\% \hbox{ CI } 14\% \hbox{ to }20\%)$}  & \cite{BCB, PL} & \pbox{15cm}{proportion of asymptomatic \\ $\times$ $(\textrm{duration of asymptomatic})^{-1}$}\\
 $\gs$ & $22\% \times  (7 \textrm{ days})^{-1}$ & \cite{OP, TBP} & \pbox{15cm}{proportion of undetected symptomatic \\ $\times$ $(\textrm{duration of symptomatic})^{-1}$} \\
 $\gh$ & \pbox{15cm}{computed by \\ observations} & \cite{TK} & discussed in Section \ref{sec:exp} \\
 $\gp$ & \pbox{15cm}{computed by \\ observations} & \cite{TK} & discussed in Section \ref{sec:exp}
\end{TAB}
\caption{Medical parameters}
\label{tbpara}
\end{table}

\subsection{The observations} \label{sec:obs}

As introduced in Section \ref{sec:int}, one of the aims of this study is to estimate the unknown parameters $\ba$ and $\bs$, and the unobservable compartments $\ia$ and $\is$. 
Because of relation \eqref{eq_on_beta}, it is clear that for the parameters only $\bs$ 
has to be evaluated.
The method that we are going to introduce in Section \ref{enkf} requires observations of some compartments in the extended SEIR model. 
In Figure \ref{pesir}, we highlight the three compartments with observations in yellow, namely $H$, $\rh$ and $\dh$. 
The data corresponding to these compartments may not be perfect, and for the analysis based 
on these data we shall take some uncertainties into consideration. 
More precisely, some random noise will be added to the observations.

For our experiment, we shall firstly and mainly use the data about Tokyo, with a total population $N=13'955'000$ (approximate mean of the population of 2020 and 2021). Subsequently, other prefectures in Japan are also considered.
As already mentioned, the compartment $H$ corresponds to the identified individuals either hospitalized or treated at home or in a hotel. 
On the other hand,  $\rh$ and $\dh$ describe the accumulated number of recovered and deceased individuals. 
The information about these three compartments are very easy to find for Tokyo, but also for any region in Japan. 
One can check for example the website of Ministry of Health, Labour and Welfare, prefecture's websites or some mass communication companies' websites to get more details. 
The information is provided on a daily basis. Note however, that these records were not very accurate at the beginning of the outbreak. 
This was caused by the delay of reports, but also by some changes in the policy for the collect of information. Our analysis will take this source of uncertainties into consideration.

\section{Ensemble Transform Kalman filter} \label{enkf}

In this section, we introduce the main tool employed in the study, namely the ensemble Kalman filter. The analysis of parameter values and unobservable compartments is performed based on it.

The method Ensemble Transform Kalman Filter (ETKF) used in this paper has been introduced
and developed  in \cite{BEM, HKS}. 
It is based on the Ensemble Kalman Filter \cite{E1, E2, E3}, but it provides analyses in a more efficient way.

Let us consider a discrete time state space model given for any time $t\in \N:=\{1,2,\ldots\}$ by
\begin{align}
\tx_{t} &= M_{t}(\tx_{t-1}) \label{eq1}\\
\ty_{t} &= H_{t}(\tx_{t}) + \epsilon_{t}, \label{obs_sys}
\end{align}
where $\tx_t$ is an $l$-dimensional state vector and $\ty_t$ an $m$-dimensional observational vector,
$M_{t}$ is the operator which defines the dynamics of the state,
$H_t$ is an operator corresponding to the observation model,
and $\epsilon_t$ provides the observation error.
All these vectors or operators are explicitly $t$ dependent.
The vector $\tx_t$ represents the state of the dynamical system, 
while $\ty_t$ is called the noisy observation, both at time $t$.
The observation error $\epsilon_{t}$ follows a normal distribution $N(\t0, \tR_{t})$, where $\tR_{t}$ is an $m \times m$ observation error covariance matrix. 
In this framework, the general question is: given a list of noisy, unbiased observations $(\ty_{t})_{t}$, how can one find the best estimates for $(\tx_{t})_{t}$\;\!?

Under the assumption of unbiased Gaussian observation error and of Gaussian initial distributions, 
one can consider the maximum likelihood approach. 
For time $t = 1, \dots, N$, the likelihood of $(\tx_{t})_{t}$ is proportional to 
\begin{align} \label{eqllh}
\prod_{t = 1}^{N} \exp\Big(-\frac{1}{2}[\ty_{t} - H_{t}(\tx_{t})]^{T}\tR_{t}^{-1}[\ty_{t} - H_{t}(\tx_{t})]\Big),
\end{align}
where the superscript $T$ denotes the transpose of a vector or of a matrix.
Clearly, maximizing \eqref{eqllh} can be transformed into a minimization problem. 
More precisely, by using equation \eqref{eq1}, one first defines the cost function
\begin{align}\label{eqco}
\sum_{t = 1}^{N}\big[\ty_{t} - H_{t}( M_{t}(\tx_{t-1}))]^{T}\tR_{t}^{-1}[\ty_{t} - H_{t}(M_{t}(\tx_{t-1}))\big].
\end{align}
Then, one ends up in looking for a list of estimates $(\tx_{t}^{a})_{t}$ of $(\tx_{t})_{t}$ which minimize \eqref{eqco}. 

When the operators $M_{t}$ and $H_{t}$ are linear, they can be represented by matrices, 
denoted also by $M_{t}$ and $H_{t}$. 
In such a situation, and with our assumption of Gaussian observation error, 
the linear model leads to the propagation of Gaussian distributions for the states.
More precisely, assume that at some time $t$, there is a prior estimate $\tx^{b}_{t}$ of the state $\tx_{t}$ with its covariance estimate $\tP_{t}^{b}$ (also called the background error covariance).
Then the Kalman filter \cite{K, KB} provides an optimal solution to minimize the variance of its uncertainty. For this reason, it is known as the minimum variance estimator. Formally, 
the prior estimate is updated by using the information of the observation $\ty_{t}$ at time $t$ as follows \cite[Eq.~(1), (2), (8) \& (10)]{HZ}\;\!:
\begin{align}
\tx_{t}^{a} &=  \tx^{b}_{t} + \tK_{t}[\ty_{t} - H_{t}\tx^{b}_{t}], \label{kaleq1}\\
\tK_{t} &= \tP_{t}^{b}H_{t}^{T}(H_{t}\tP_{t}^{b}H_{t}^{T} + \tR_{t})^{-1},\label{kaleq2}\\
\tP_{t}^{a} &= (I - \tK_{t} H_{t})\tP_{t}^{b}, \label{kaleq3}\\
\tP_{t+1}^{b} &=M_{t+1}\tP_{t}^{a}M_{t+1}^{T}, \label{kaleq4} 
\end{align}
where $\tK_{t}$ is called the Kalman gain, $\tx_{t}^{a}$ the analysis, and $\tP_{t}^{a}$ its estimated covariance. 
With the analysis, one can use the forecast model \eqref{eq1} to evolve the analysis and get the prior estimate and the covariance at the next time step.

When the operators $M_{t}$ and $H_{t}$ are nonlinear, the Kalman filter is not applicable
in its original form, but a simple extension with linear approximation is known to be effective. The straightforward nonlinear extension of the Kalman filter is known as the extended Kalman filter (EKF). 
In this case, $M_{t+1}$ and $H_{t}$ in equation \eqref{kaleq1} to \eqref{kaleq4} are replaced by their linear approximations respectively. 
For low dimensional models, the cost of EKF is low. However, when the model complexity and the state dimension increases, one may consider the ensemble Kalman filter. 
In this case, an ensemble of model states is used to estimate the covariance matrices.
Namely, let $\{\tx_{t-1}^{a(i)}, i = 1, \dots, k\}$ be an ensemble of estimates of the model states at time $t - 1$. By letting this ensemble evolve according to the model \eqref{eq1}, one gets the forecast or background ensemble $\{\tx_{t}^{b(i)}, i = 1, \dots, k\}$ at time $t$, where
\begin{align*}
\tx_{t}^{b(i)} = M_{t}(\tx_{t-1}^{a(i)}).
\end{align*}
The mean of the ensemble can be used as the most probable estimate of the state, with its covariance estimated by the sample covariance of the ensemble:
\begin{equation*}
\bar{\tx}_{t}^{b} = \frac{1}{k}\sum_{i = 1}^{k}\tx_{t}^{b(i)}
\qquad \hbox{ and }\qquad
\tP_{t}^{b} =  \frac{1}{k-1}(\tX_{t}^{b})(\tX_{t}^{b})^{T},
\end{equation*}
where $\tX_{t}^{b}$ is an $l \times k$ matrix with the $i^{\rm th}$ column $\tx_{t}^{b(i)} - \bar{\tx}_{t}^{b}$.

The next step is to update the forecast ensemble mean and its covariance by using the information of the observation $\ty_{t}$. 
The updated ensemble mean is denoted by $\bar{\tx}_{t}^{a}$ and its covariance by $\tP^{a}_t$. 
To get them, one first needs to determine an analysis ensemble $\{\tx_{t}^{a(i)}, i = 1, \dots, k\}$.
Similar to the background ensemble, the mean and covariance of the analysis ensemble are given by 
\begin{align}
\bar{\tx}_{t}^{a} &= \frac{1}{k}\sum_{i = 1}^{k}\tx_{t}^{a(i)}, \label{eqam}\\
\tP_{t}^{a} &=  \frac{1}{k-1}(\tX_{t}^{a})(\tX_{t}^{a})^{T}. \label{eqap}
\end{align}
In our application the operator $H_t$ is linear, so that $H_t$ is given by a matrix.
In this case and for ETKF, the computations can be written as follows \cite[Eq.~(15)$\sim$(17)]{HZ}\;\!:
\begin{align}
\bar{\tx}_{t}^{a} &= \bar{\tx}_{t}^{b} + \tX_{t}^{b}\tilde{\tP}_{t}^{a}(H_{t}\tX_{t}^{b})^{T}\tR_{t}^{-1}\big(\ty_{t} -\overline{H_{t}\tx^{b}_{t}}\big),\\
\tilde{\tP}_{t}^{a} &= \big[(k - 1)I + \big(H_{t}\tX_{t}^{b}\big)^{T}\tR_{t}^{-1}H_{t}\tX_{t}^{b}\big]^{-1},\\
\tX_{t}^{a} &= \tX_{t}^{b}\big[(k -1)\tilde{\tP}_{t}^{a}\big]^{\frac{1}{2}}\label{eqwei},
\end{align}
where the $k\times k$ matrix $\tilde{\tP}_{t}^{a}$ is called the analysis error covariance in ensemble space. For \eqref{eqwei}, the background perturbation 
$\tX_{t}^{b}$ is transformed into the analysis perturbation $\tX_{t}^{a}$ with the weight $[(k -1)\tilde{\tP}_{t}^{a}]^{\frac{1}{2}}$.
By adding $\bar{\tx}_{t}^{a}$ to each column of $\tX_{t}^{a}$, one gets the analysis ensemble $\{\tx_{t}^{a(i)}, i = 1, \dots, k\}$ that satisfies \eqref{eqam}. 
Note that the ETKF is a deterministic filter since no randomly perturbed observations are used in the computation \cite{WH}, and it is also a square-root filter because its takes the power $\frac{1}{2}$ of the matrix $\tilde{\tP}_{t}^{a}$ \cite{TAB}.

To avoid a variance underestimation,  caused for example by the limited ensemble size,
an artificial inflation of the ensemble spread is usually applied. 
In \cite{HZ}, several ways to perform the variance inflation are discussed. 
In our study, we apply multiplicative inflation and additive inflation. 
For multiplicative inflation, the background error covariance $\tP_{t}^{b}$ is multiplied by a tunable factor $\rho > 1$ before the analysis.
The multiplicative inflation can be thought as a procedure to increase the influence of the current observations on the analysis. 
For additive inflation, since the extended SEIR model is a conserved closed system, 
it is then applied by adding a random vector to the background ensemble. 
The random vector can be sampled from the normal distribution with mean \textbf{0} and some covariance matrix. 
This covariance matrix is assumed to be proportional to the background error covariance $\tP_{t}^{b}$ by a tunable parameter $\alpha \in (0, 1)$ \cite{HZ}.

For the current study with a low dimensional model, we could use EKF or EnKF. However, by considering the possibility of increasing the model complexity in future, we opted for ETKF with 50 ensemble members.

\section{Experiments} \label{sec:exp}

In this section, we explain the design of the experiments and illustrate the outcomes
by using observations from Tokyo Metropolis.
These observations ($H$, $R$, and $D$) are obtained from the website 
\cite{TK}.
Discussions and comparisons with other prefectures are provided 
in the subsequent section.
As mentioned in Section \ref{sec:obs}, health officials started providing data from February 2020 but they included some uncertainties caused by delay or by policy changes. Our experiments start with the observations from March $6^{\mathrm{th}}$, 2020, when the record of data became more systematic.

\subsection{The experiments, using data from Tokyo}\label{sec:tok}

For the state space model mentioned in Section \ref{enkf}, the dynamical system \eqref{eq1} is described by the differential system \eqref{esirm}. 
To estimate the parameter $\bs$, we use the augmented state by adding one more equation $\frac{d\bs}{dt} = 0$ to the system \eqref{esirm}, assuming persistence for $\bs$. 
Namely, $\bs$ will stay at the same value during the time integration process, and will be updated during the analysis step of the data assimilation.
In other words, if $\bs$ has a correlation with the observations, $\bs$ will be updated together with the other states.
In system \eqref{esirm}, the unit of time $t=1$ represents one day. 
Thus, the one day forecast from day $n$ is obtained by integrating the system \eqref{esirm} on $[n, n+1]$ with the initial values equal to the analysis on day $n$.
To avoid negative values in the analysis step, all the compartments are transformed to $\log$ scale with base $\e$.
The lower case will be used for these new variables, as for example $s(t) := \log S(t)$, 
and the corresponding equation becomes 
$\frac{ds(t)}{dt} = \frac{1}{S(t)} \frac{dS(t)}{dt}$. 
The DA analysis update directly applies to the $\log$-transformed values.

Now the forecast value of $\tx_{t}$ is a $9$-dimensional vector 
\begin{equation}\label{eq_forecast_vector}
\big(e(t), \iia(t), \iis(t), h(t), \rrh(t), \ddh(t), \rra(t), \rrs(t), \log \bs\big)^{T}.
\end{equation}
Note that the compartment $s(t)$ is not explicitly included since it can be
deduced from the conservation relation \eqref{eq_conservation}.
Before performing the DA analysis update, the above vector has to be projected to the observation space by the observation system \eqref{obs_sys}.
In our setting, the operator $H_t$ introduced in \eqref{obs_sys} is defined by the 
projection which sends the forecast vector \eqref{eq_forecast_vector} on 
$\big(h(t), \rrh(t), \ddh(t)\big)^{T}$, namely on the three compartments with observations. 

To initiate the experiments, initial values for all compartments and parameters have to be provided. 
When these initial values are far away from the observations, the system goes through an unreasonable transition period. To avoid or reduce the unreasonable transition, preliminary experiments were performed for determining suitable initial values.
Namely, we run the system for a few days, starting with a few individuals in some compartments and a presumed value for $\bs$.
Different values of $\bs$ were tested, until values comparable to the data of the compartments $H, \dh, \rh$ on day zero (March $6^{\mathrm{th}}$, 2020) were obtained. The corresponding values for the different compartments were then chosen as initial conditions.
However, independent normal distributed random errors are also added to create 50 ensemble members 
$\{\tx^{a(i)}_0, i = 1, \dots, 50\}$.
Note that each member $\tx^{a(i)}_0$ is a $9$-dimensional vector containing the initial conditions for the compartments and for the parameter $\bs$. 

The ensemble on day zero are then integrated by the model \eqref{esirm} for one day and the ETKF will update the one day forecast by using the observations. The procedure of the experiments is shown in Figure \ref{pesir2}.

\begin{center}
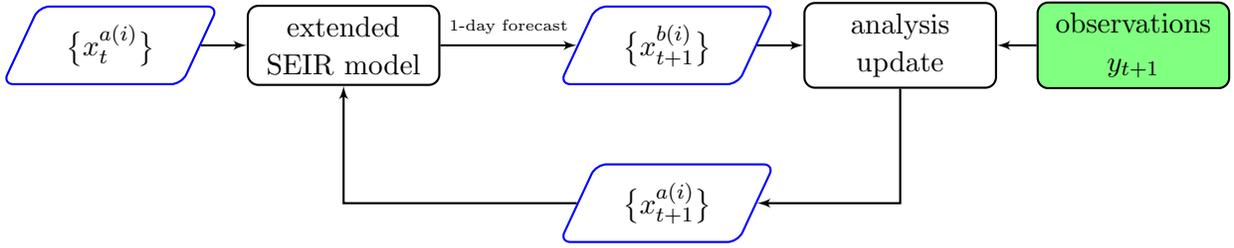

\begin{tikzpicture}[node distance=3.1cm, auto]
\node [block1, draw = blue] (init) {\small $\big\{\tx^{a(i)}_t\big\}$};
\node [block2, draw = black, right of=init] (E) {\small extended \\ SEIR model};
\node [block1, draw = blue, right of=E, xshift=1.2cm] (I1) {\small $\big\{\tx^{b(i)}_{t+1}\big\}$};
\node [block2, draw = black, right of=I1] (I2) {\small analysis update};
\node [block2, fill=green!50, right of= I2] (H) {\small observations $\ty_{t+1}$};
\node [block1, draw = blue, below of = I1, yshift=1cm] (R2) {\small $\big\{\tx^{a(i)}_{t+1}\big\}$};
\path [line] (init) -- node {}(E);
\path [line] (E) -- node {\tiny 1-day forecast}(I1);
\path [line] (I1) -- node {}(I2);
\path [line] (H) -- node {}(I2);
\path [line] (I2) |- node [near end, xshift = -0.8cm] {}(R2);
\path [line] (R2) -| node [near end, xshift = -0.8cm] {}(E);
\end{tikzpicture}
\captionof{figure}{Data assimilation flow-chart} 
\label{pesir2}
\end{center}

\begin{figure}[hbtp]
\centering
\includegraphics[width=16cm, height =6cm]{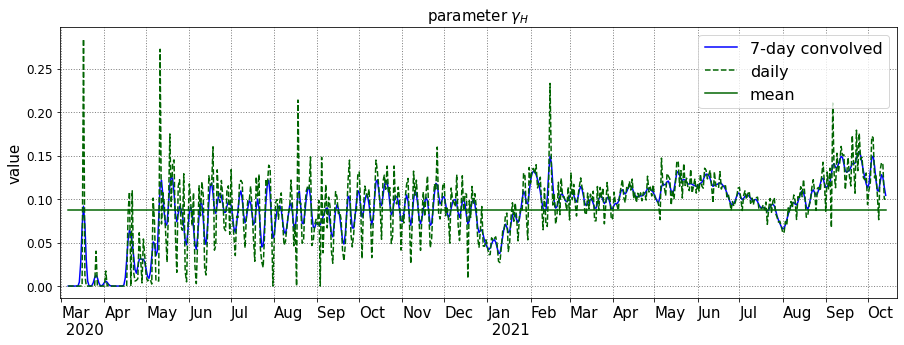}
\includegraphics[width=16cm, height =6cm]{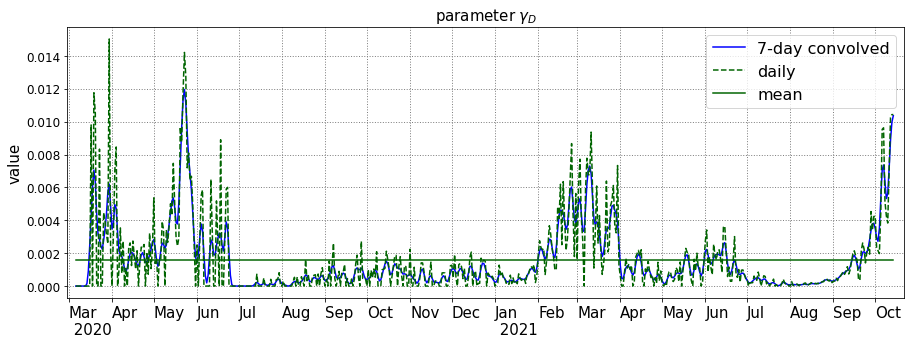}
\caption{Daily estimation of the parameters $\gh$ and $\gp$}
\label{fig:gg}
\end{figure}

As mentioned in Section \ref{sec:mpara}, parameters $\gh$ and $\gp$ are estimated by assimilating the observations. 
Consider equation $\frac{d \rh}{dt} =  \gh H$ in model \eqref{esirm}. Since the time is discrete, we can rewrite this relation as $\rh(t+1)-\rh(t) = \gh(t)H(t)$ which leads to
\begin{align}\label{eq:gh}
\gh(t) = \frac{\rh(t+1) - \rh(t)}{H(t)}.
\end{align}
Similar computation can be generated for $\gp$ also, namely, 
\begin{align}\label{eq:gp}
\gp(t) = \frac{\dh(t+1) - \dh(t)}{H(t)}.
\end{align}
The daily values of both parameters are shown in Figure \ref{fig:gg}.
For the implementation of these parameters for the computation of the $1$-day forecast, 
we have used a slightly smoothed version obtained by a $7$-day convolution with the symmetric weights $\frac{1}{64}\big(1,6,15,20,15,6,1)$. 
The effect is to decrease the amplitude of the weekly oscillations. 

For the integration process, we have also included some uncertainties to all parameters: to the ones presented in Table \ref{tbpara}, but also to the values of the parameters
$\gh$, $\gp$.
Thus, the $1$-day forecast is obtained with the parameters of the previous day, each of them perturbed by a normal distribution $N\big(0, (M/10)^2\big)$, 
where $M$ corresponds to the value of this parameter.
These perturbations are independent and randomly generated for each $1$-day forecast.

The last initial setting for data assimilation is the observation error covariance matrix $\tR_{t}$. 
We assume that the covariance of the errors between different observations is zero, 
namely, the matrix $\tR_{t}$ is diagonal.
Since we observe three states, $\tR_{t}$ is a $3 \times 3$ diagonal matrix.
Clearly, the bigger the variance, the weaker contribution the corresponding observation will make to the update.

For simplicity, we assume that the observation error covariance is independent of time. 
The diagonal elements are chosen as $\big(\log(1.3)\big)^{2}$ for any time $t$, representing 
the observation error variance: Namely, for $x\in \{h, \rrh,\ddh\}$ one has $\big[x-\log(1.3), x+\log(1.3)\big]$ as $68\%$ CI and $\big[x-2\log(1.3), x+2\log(1.3)\big]$ as $95\%$ CI.
Considering that all the observations have been transformed into $\log$ scale, that is,
the observations in original scale are distributed in $\big[X/(1.3) , (1.3)X\big]$ as $68\%$ CI and $\big[X/(1.69), (1.69) X\big]$ as $95\%$ CI.

\begin{figure}[hbtp]
\centering
\includegraphics[width=5cm]{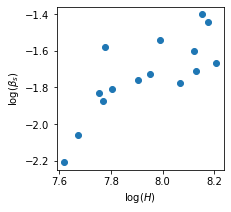}
\caption{Scatter plot of $h$ and $\log(\bs)$ on November 17, 2020}
\label{fig:sca}
\end{figure}

Before showing the analysis result, a scatter plot of the values of $h$ and of $\log(\bs)$ for the different ensemble members on a certain day is provided in Figure \ref{fig:sca}.
Indeed, for the parameter estimation, one assumption is the linear relation between the ensemble of parameters and the ensemble of compartments with observations. 
As shown in Figure \ref{fig:sca}, a positive correlation is observed between $h$ and $\log(\bs)$. Thus we can use the observation of $H$ for the estimation of $\bs$.

\begin{figure}[hbtp]
\centering
\includegraphics[width=11cm]{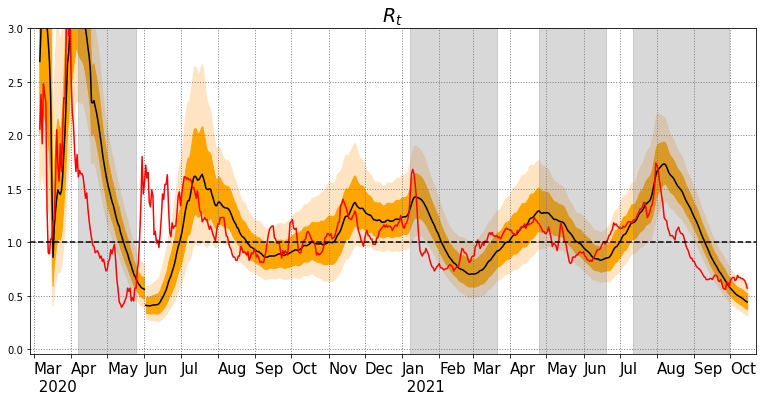}
\caption{Analysis value for $\rt$, and comparison with the value of Toyokeizai.net }
\label{fig:br}
\end{figure}

In Figure \ref{fig:br} we finally provide the analysis value for $\rt$ from March $6^{\mathrm{th}}$ 2020 to the middle of October 2021.
The black curve represents the mean value of the ensemble. 
The two shaded regions represent the $68\%$ CI (dark orange) and $95\%$ CI (light orange). The analysis value of $\rt$ is computed by using $\bs$ together with the information contained in \eqref{eqro} and in \eqref{eq_on_beta}. The red curve shows the $\rt$ provided by Toyokeizai.net for reference. 
The grey shaded regions represent the periods of state of emergency in Tokyo.
Note that the jump at the end of May 2020 is due to an abrupt change in the value of $\tau_H$ as indicated in Table \ref{tbpara}.
Additional discussions about this graph will be provided in the subsequent section.

\begin{figure}[hbtp]
\centering
\subfigure{\includegraphics[width=7.8cm]{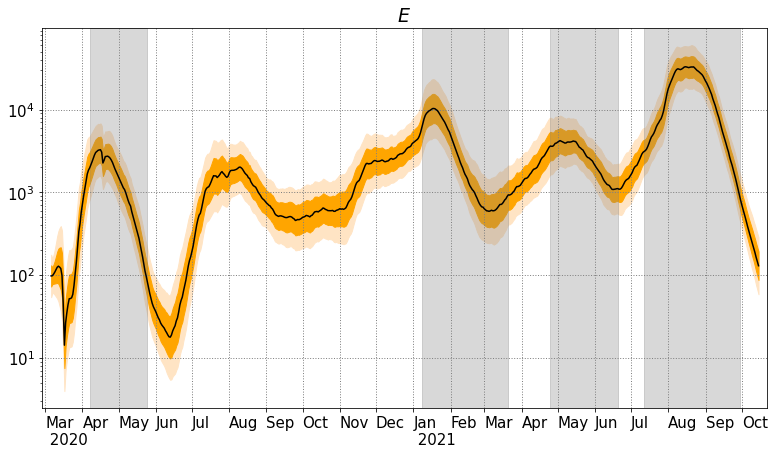}}
\subfigure{\includegraphics[width=7.8cm]{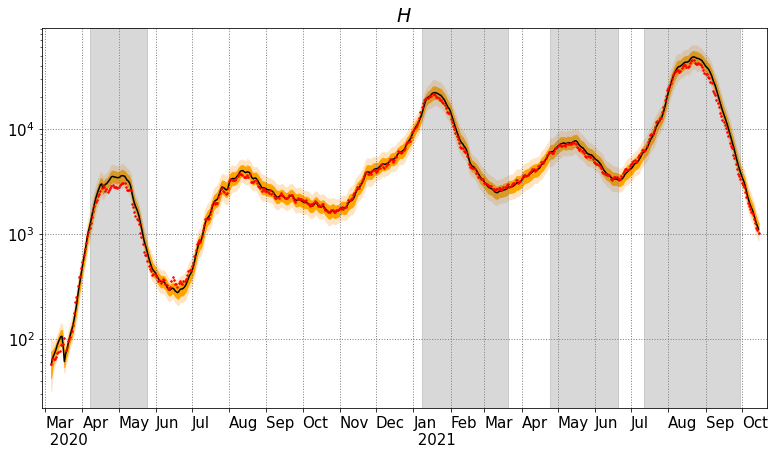}}
\subfigure{\includegraphics[width=7.8cm]{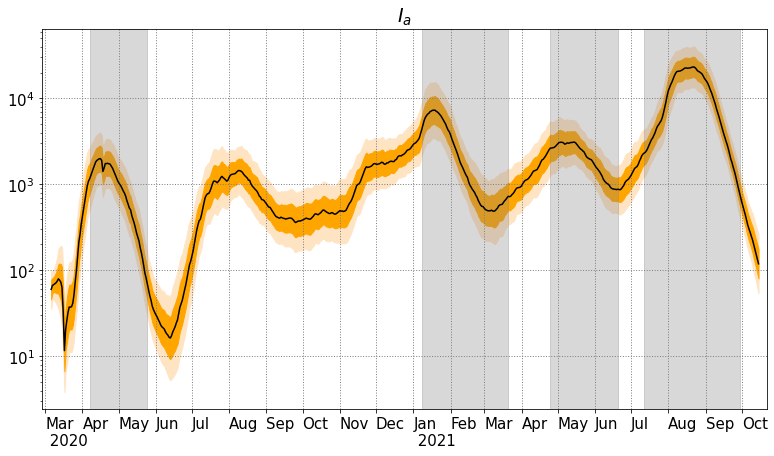}}
\subfigure{\includegraphics[width=7.8cm]{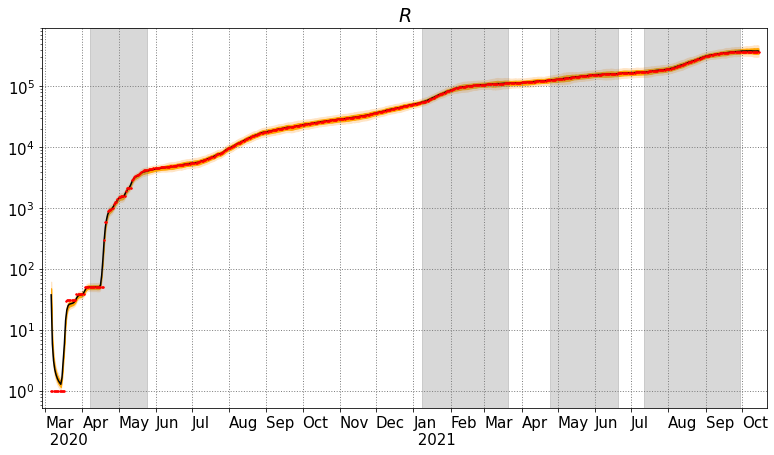}}
\subfigure{\includegraphics[width=7.8cm]{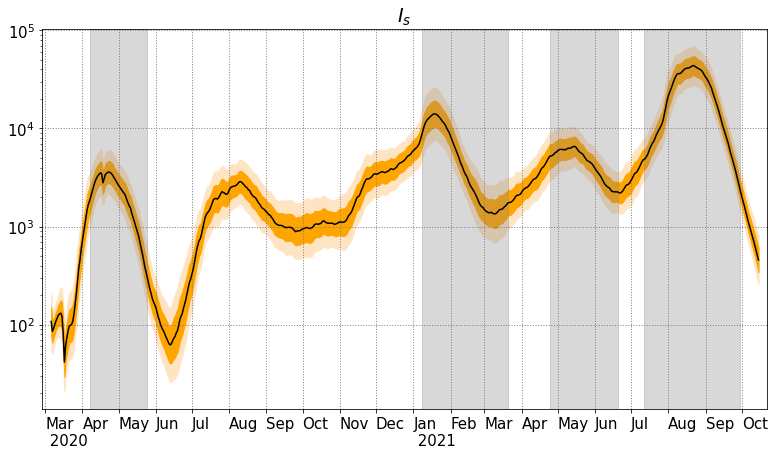}}
\subfigure{\includegraphics[width=7.8cm]{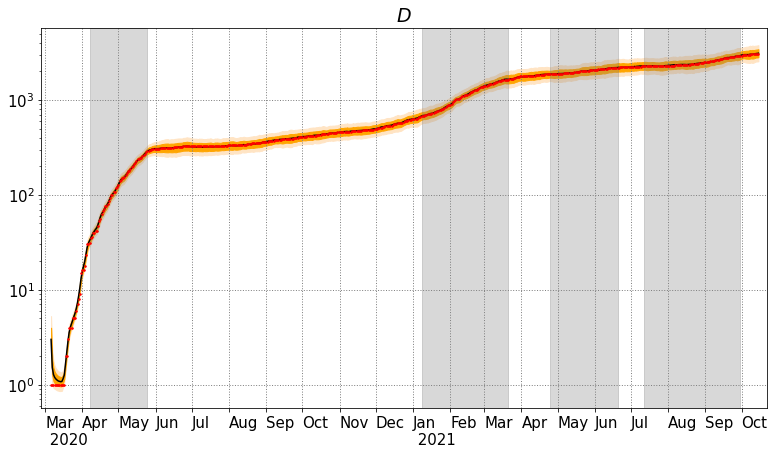}}
\caption{The 6 compartments for Tokyo}
\label{fig:to6}
\end{figure}

The analysis results for the compartments $E$, $\ia$, $\is$, $H$, $\rh$, $\dh$ with confidence intervals are also presented in Figure \ref{fig:to6}. For compartments $H$, $\rh$, $\dh$, the red dots represent the observations of Tokyo. 
For all pictures, the $y$-axis is in $\log_{10}$ scale.

\subsection{Technical discussion, using data from Tokyo}

In Figure \ref{fig:br}, one first notices that the similarity between the analysis $\rt$ and the reference  $\rt$ provided by Toyokeizai.net. 
The Toyokeizai $\rt$ on day $t$ is computed by the following formula:
\begin{align*}
\bigg(\frac{\sum_{i = t - 6}^{t} \hbox{new confirmed cases on day $i$}}{\sum_{i = t - 13}^{t - 7} \hbox{new confirmed cases on day $i$}}\bigg) ^ {\frac{\hbox{average generation time}}{\hbox{reporting interval}}},
\end{align*}
where the average generation time is assumed to be 5 days and the reporting interval is assumed to be 7 days. 
The formula is a simplified version of a maximum likelihood estimation for the effective reproduction number provided by H. Nishiura \cite{NH}. 
Note that this formula is based on the daily confirmed cases, which experiences quick changes everyday. 
On the other hand, the analysis $\rt$ and its confidence interval is computed by formula \eqref{eqro} which involve the analysis ensemble of $\bs$ and $\ba$. 
The average over the ensemble member provides a smoother evolution of $\rt$, which is certainly closer to a true statistic's evolution.
Note that a third approach for the computation of $\rt$ is available in 
\cite{SRM}: it involves an agent-based model together with a particle filter.

In Figure \ref{fig:to6}, we show the analysis mean of compartments $E$, $\ia$, $\is$ with the $68\%$ CI (dark orange) and $95\%$ CI (light orange), 
and the same results for compartments $H$, $\rh$ and $\dh$ with the daily observations represented by red dots. 
One notes that the analyses of $H$, $\rh$, $\dh$ fit the observations appropriately except at the beginning of the observation period.
This may be related to the initial guess and the model imperfections. 

\begin{figure}[hbtp]
\centering
\subfigure{\includegraphics[width=9.2cm]{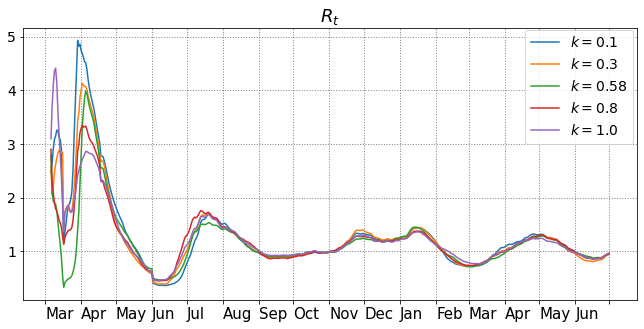}}
\caption{Comparison of analysis results with different ratio $k$ between $\ba$ and $\bs$}
\label{fig:comp}
\end{figure}

\begin{figure}[hbtp]
\centering
\subfigure{\includegraphics[width=16cm]{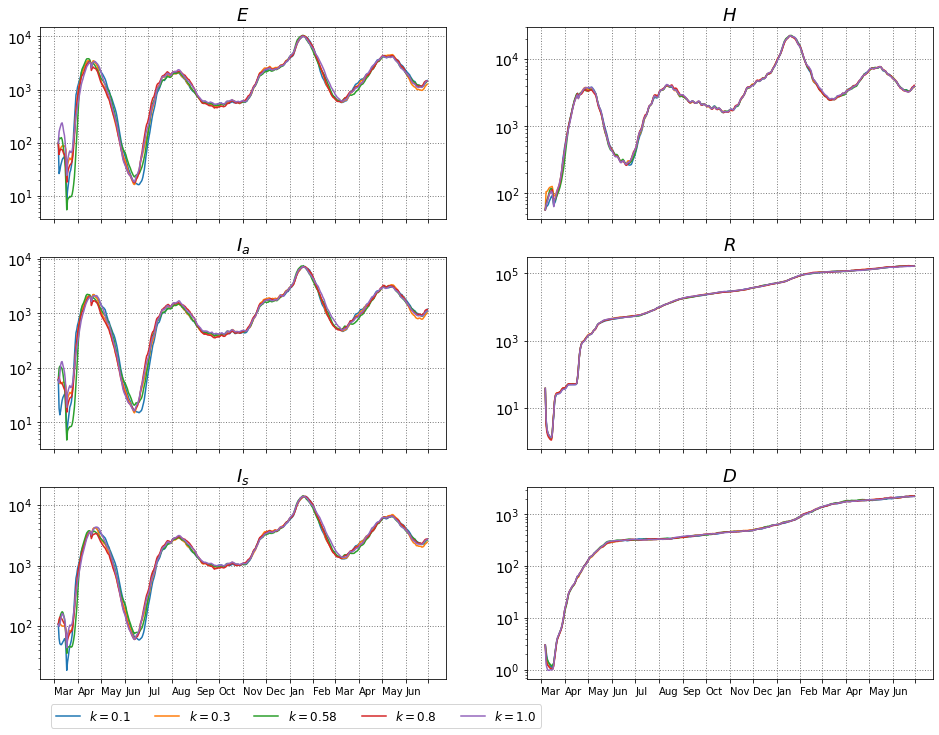}}
\caption{Analysis mean of compartments with different ratio $k$}
\label{fig:otherk}
\end{figure}

As mentioned in Section \ref{sec:mpara}, the ratio between $\ba$ and $\bs$ is fixed at $0.58$ based on the information provided by the literature.
Since this ratio is quite uncertain, we performed similar investigations with 
different ratios varying from $0.1$ to $1.0$.
For $k=\frac{\ba}{\bs}$ with $k = 0.1, 0.3, 0.58, 0.8$ and $1.0$, we run the experiments independently (always with the data from Tokyo) and the analysis mean of $\rt$ are shown in Figure \ref{fig:comp}. 
The patterns of each $\rt$ are similar, especially after May $31^{\mathrm{st}}$, 2020. 
However, at the beginning of the outbreak and during the first state of emergency, all the analysis value of $\rt$ experience quick changes but with different speed.
In Figure \ref{fig:otherk} we also provide the analysis mean of the different compartments, since they also depend on the ratio between $\ba$ and $\bs$. The same conclusion is obtained: except for the first three months, the possible ratio do not generate any noticeable difference. 

\begin{figure}[hbtp]
\centering
\subfigure[]{\includegraphics[width=8.2cm]{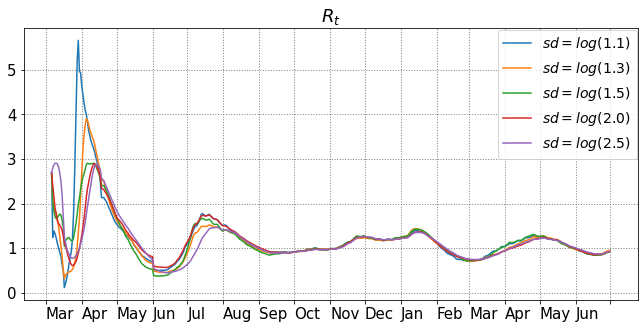}}
\subfigure[]{\includegraphics[width=8.2cm]{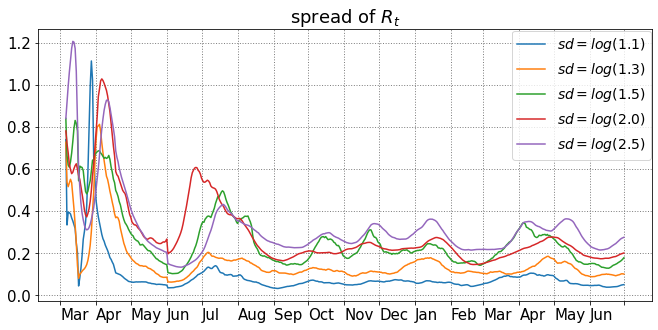}}
\caption{Analysis mean and spread of $\rt$ with different observation covariance}
\label{fig:rt_ms}
\end{figure}

\begin{figure}[hbtp]
\centering
\subfigure[]{\includegraphics[width=14cm,height = 9.5cm]{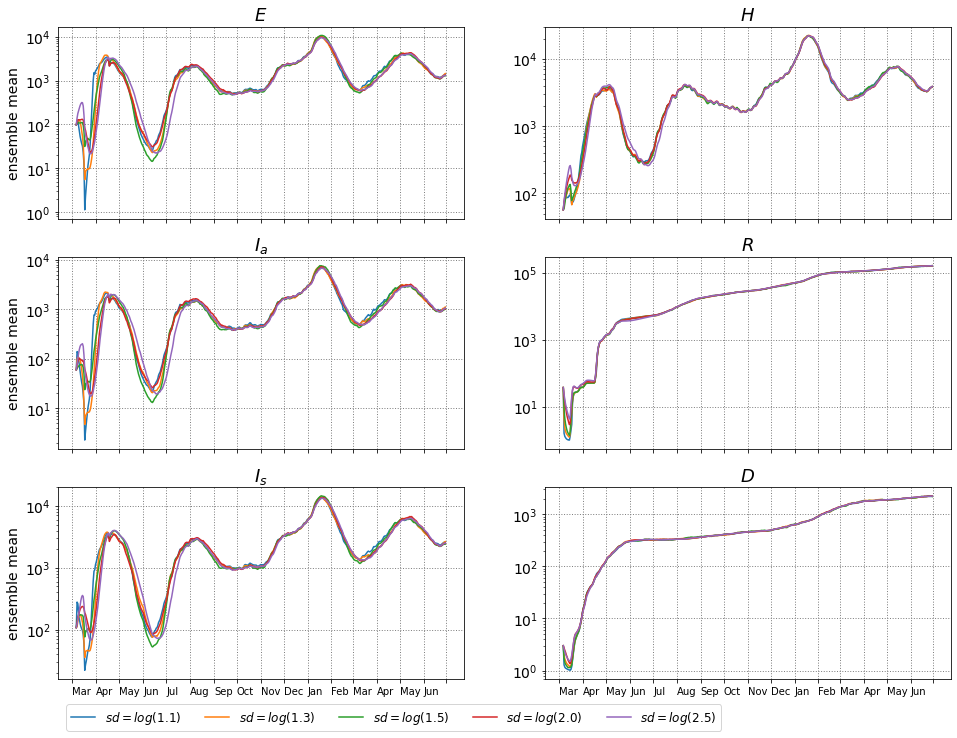}}
\subfigure[]{\includegraphics[width=14cm,height = 9.5cm]{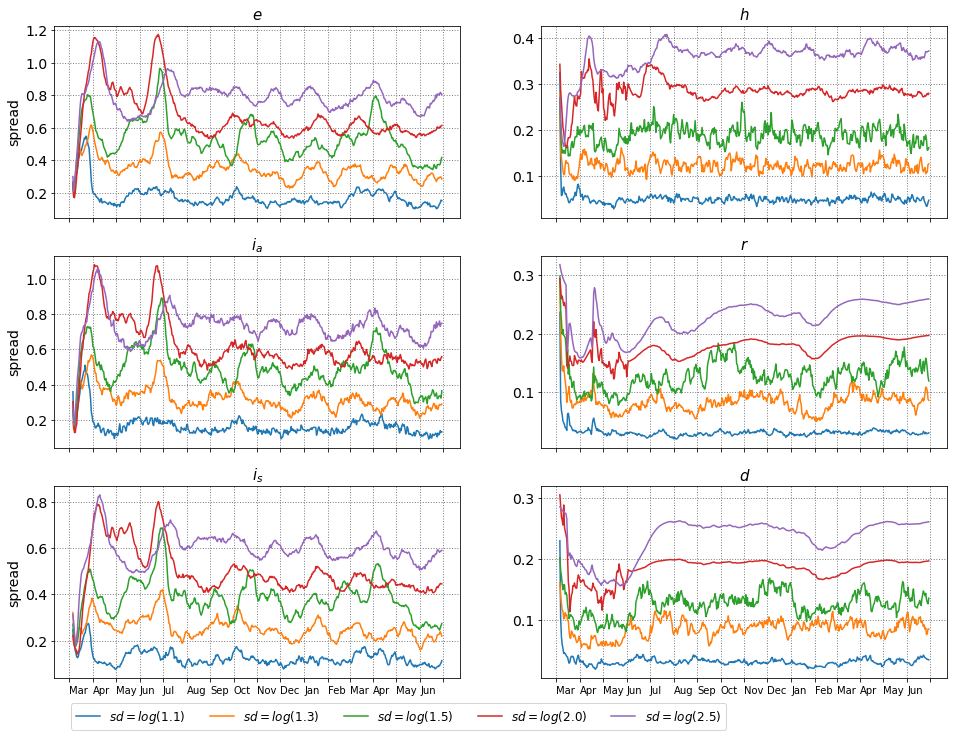}}
\caption{Analysis mean (a) and spread (b) of the compartments with different observation covariance}
\label{fig:oth_ms}
\end{figure}

In Section \ref{sec:tok}, the error covariance matrix  $\tR_{t}$ is defined as a diagonal matrix $\big((\log(1.3)\big)^2 I$ where $I$ is the identity  $3 \times 3$ matrix. 
We also run a sensitivity test for different choices of the covariance matrix, and checked if it has an effect on the analysis results. 
We choose $\tR_{t}$ equals to ${\rm sd}^2 I$ with ${\rm sd} = \log(1.1)$, $\log(1.3)$, $\log(1.5)$, $\log(2.0)$, and $\log(2.5)$. For these investigations 
the ratio between $\ba$ and $\bs$ is fixed as $0.58$. In Figure \ref{fig:rt_ms}(a), 
the analysis means are provided, and they are clearly close to each other except at the beginning stage of the outbreak. The different spreads are shown in 
Figure \ref{fig:rt_ms}(b), 
and one notices that the spread increases with the ${\rm sd}$, but for ${\rm sd}$ larger than $\log(1.5)$, the increase is less clear.
We also show the analysis mean of compartments $E$, $\ia$, $\is$, $H$, $\rh$ and $\dh$ in Figure \ref{fig:oth_ms}(a). 
Again, the analysis results of these compartments with different observation error covariance are quite close to each other. 
In Figure \ref{fig:oth_ms}(b), the spread is computed by the standard deviation of ensemble members of the $\log-$transformed compartments. 
The relations of spread are quite clear for the three compartments with observations, 
and there are some overlaps for spread of $e$, $i_a$ and $i_s$, but in general the spread increases with ${\rm sd}$.

\begin{figure}[hbtp]
\centering
\subfigure[]{\includegraphics[width=7.6cm]{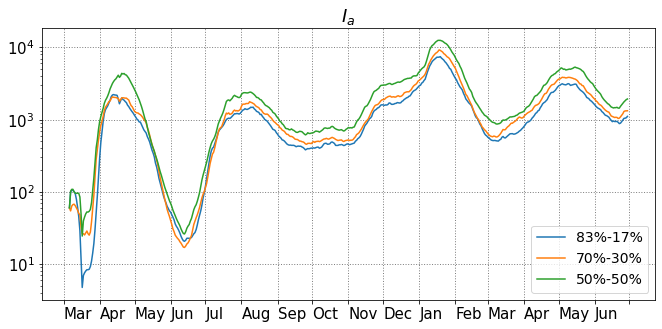}}
\subfigure[]{\includegraphics[width=7.6cm]{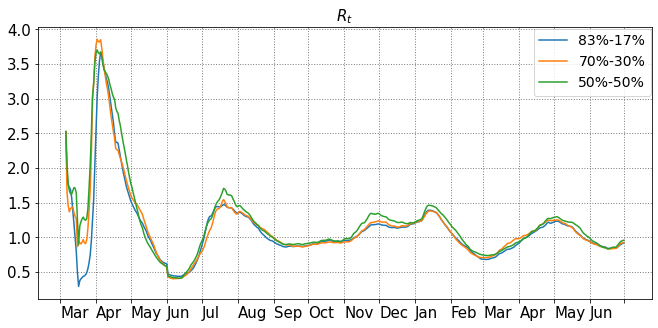}}
\caption{Analysis mean of state $\ia$ (a) and $\rt$ (b) with different ratios between symptomatic and asymptomatic}
\label{fig:asy_portion}
\end{figure}

As shown in Table \ref{tbpara}, the proportion of symptomatic and asymptomatic used in our investigation is $83\%$ and $17\%$, respectively. 
These values were determined based on sources mentioned in Section \ref{sec:mpara}.
However, since asymptomatic cases are very difficult to detect, and since we can not be fully confident in the above ratio, a sensitivity test is necessary. To do this, we choose two different combinations which are $70\%$ and $30\%$, and 
$50\%$ and $50\%$, and keep the observation error covariance matrix at $\big((\log(1.3)\big)^2 I$ and $k = 0.58$. 
Compared to the previous settings, in these two scenarios more people become asymptomatic and recover without showing any symptoms, as illustrated in Figure \ref{fig:asy_portion}(a). 
On the other hand, the three curves for $\rt$ are close to each others, 
but the one corresponding to only $50\%$ of symptomatic is slightly higher than the other two, see Figure \ref{fig:asy_portion}(b).
In order to understand this, let us recall that the relation between the transmission coefficients $\ba$ and $\bs$ is set to be $0.58$, it means that the asymptomatic
cases are less infectious than the symptomatic and symptomatic cases. 
Thus, when more people do not show any symptoms, 
the transmission coefficient has to be bigger to create enough cases to fit the observations. As a consequence, $\rt$ will also be bigger.

\section{Discussion, comparisons, and future projects}\label{sec:disc}

Let us start by a few easy observations about the results obtained in the previous section.

The analysis values of $\rt$ quickly decrease during the first state of emergency, 
and this parameter is lower than $1$ around the middle of May 2020. This indicates 
that the disease has started to die out. 
During the second state of emergency, analysis values of $\rt$ had a much slower decay, 
but nevertheless in early February 2021, $\rt$ successfully drops below $1$. 
However it starts to increase from March while the state of emergency had not been released. 
In the third state of emergency, the analysis $\rt$ had the slowest decay, compared
to the previous two states of emergency. 
During the first several weeks of the fourth state of emergency, the
effective reproduction number continued increasing, even at a faster speed.
It is only after the end of the first week of August that a decay finally took place.

For $E$, $\ia$ and $\is$, one finds that the three values increase at the beginning of the first state of emergency, then they slow down and change the direction to decrease. 
Similar behaviors can be observed in the second state of emergency, 
and the values start to increase around the middle of March 2021 while the $\rt$ starts yo increase at the beginning of March. There is clearly a delay between the changes of $\rt$ and their effect on the three compartments $E$, $\ia$ and $\is$.

Let us now provide a comparison with two different regions from Japan, namely Osaka and Kanagawa, and also provide the analysis for Japan as a whole.
The experiments for Osaka, Kanagawa, and Japan start from March $26^{\mathrm{th}}$ 2020, March $18^{\mathrm{th}}$ 2020, March $1^{\mathrm{st}}$ 2020 respectively. 
The observed compartments of these regions are the same as those for Tokyo, namely $H$, $\dh$ and $\rh$.
The way to determine the initial setting is the same as that for Tokyo and we use the same setting for the parameters which have been presented in Table \ref{tbpara}. 
The daily values for the parameter $\gh$ and $\gp$ for each region are also computed by using \eqref{eq:gh} and \eqref{eq:gp}. 
The ratio between $\ba$ and $\bs$ is $0.58$, and the observation error covariance matrix is assumed to be $\big(\log(1.3)\big)^2 I$.

\begin{figure}[hbtp]
\centering
\subfigure[Osaka]{\includegraphics[width=8cm,height = 5cm]{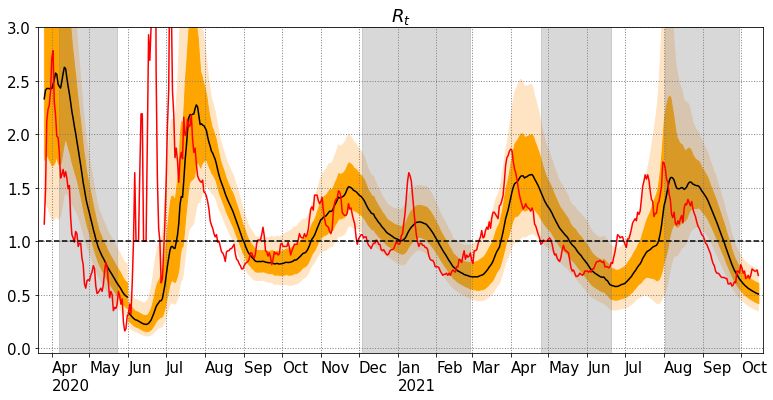}}
\subfigure[Kanagawa]{\includegraphics[width=8cm,height = 5cm]{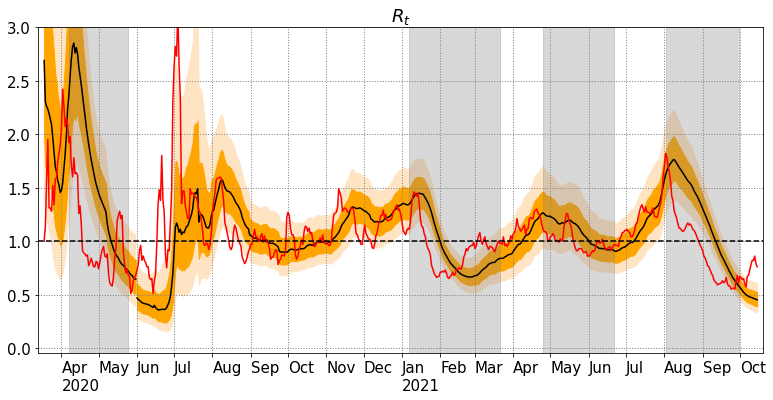}}
\subfigure[Japan]{\includegraphics[width=8cm,height = 5cm]{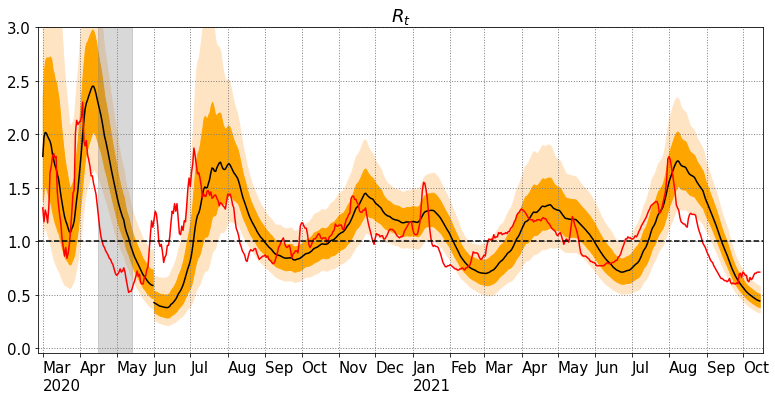}}
\caption{Analysis value of $\rt$ for Osaka, Kanagawa and Japan}
\label{fig:br_o}
\end{figure}

\begin{figure}[hbtp]
\centering
\subfigure{\includegraphics[width=17cm,height = 9.5cm]{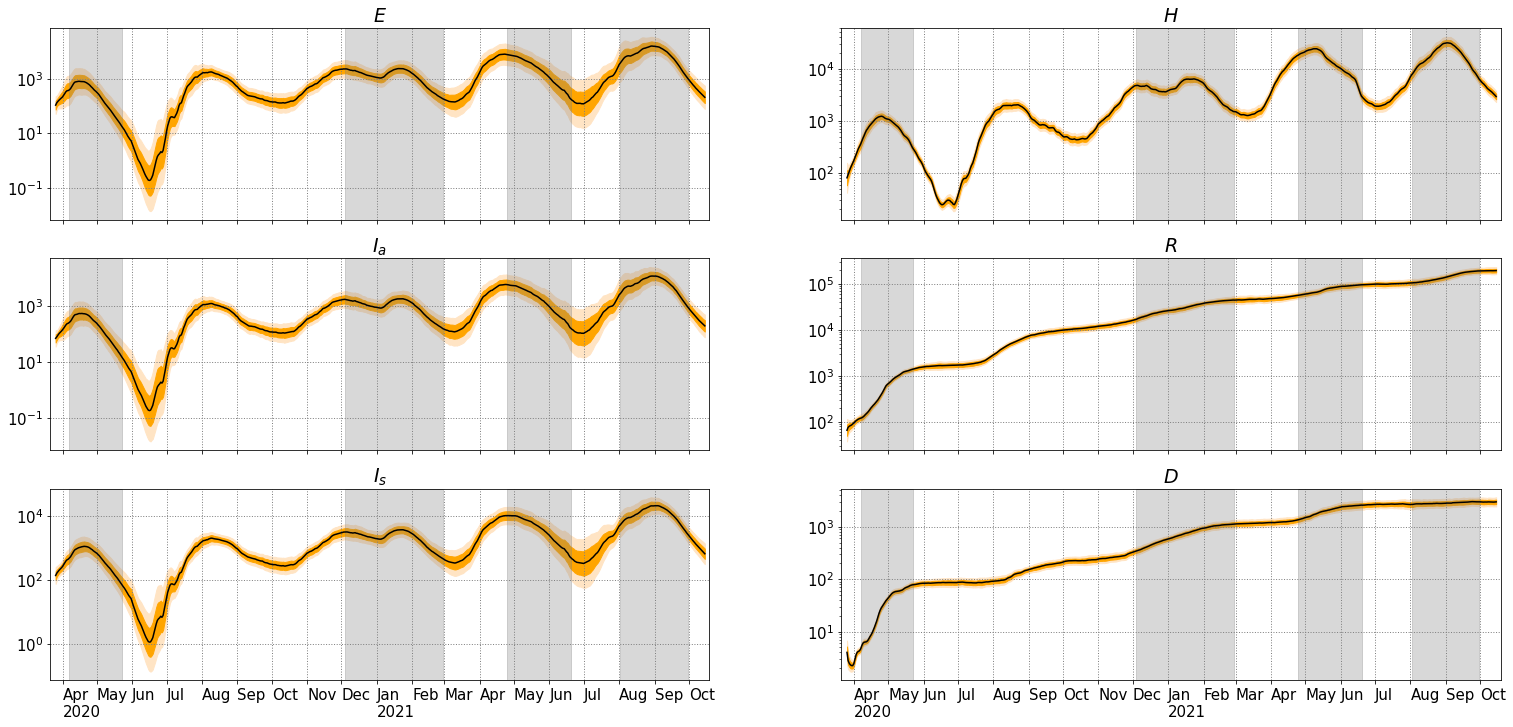}}
\caption{Analysis results of Osaka}
\label{fig:com_os}
\end{figure}

\begin{figure}[hbtp]
\centering
\subfigure{\includegraphics[width=17cm,height = 9.5cm]{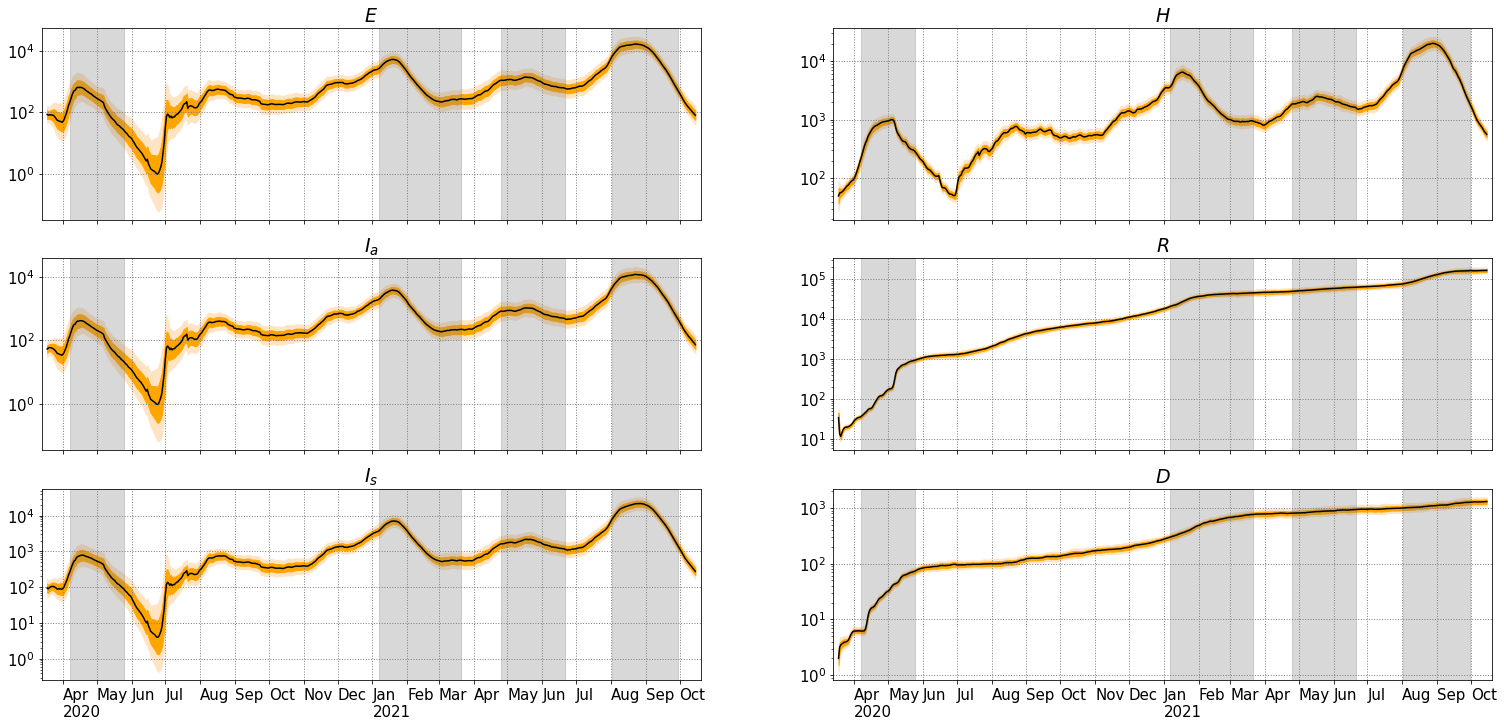}}
\caption{Analysis results of Kanagawa}
\label{fig:com_ka}
\end{figure}

\begin{figure}[hbtp]
\centering
\subfigure{\includegraphics[width=17cm,height = 9.5cm]{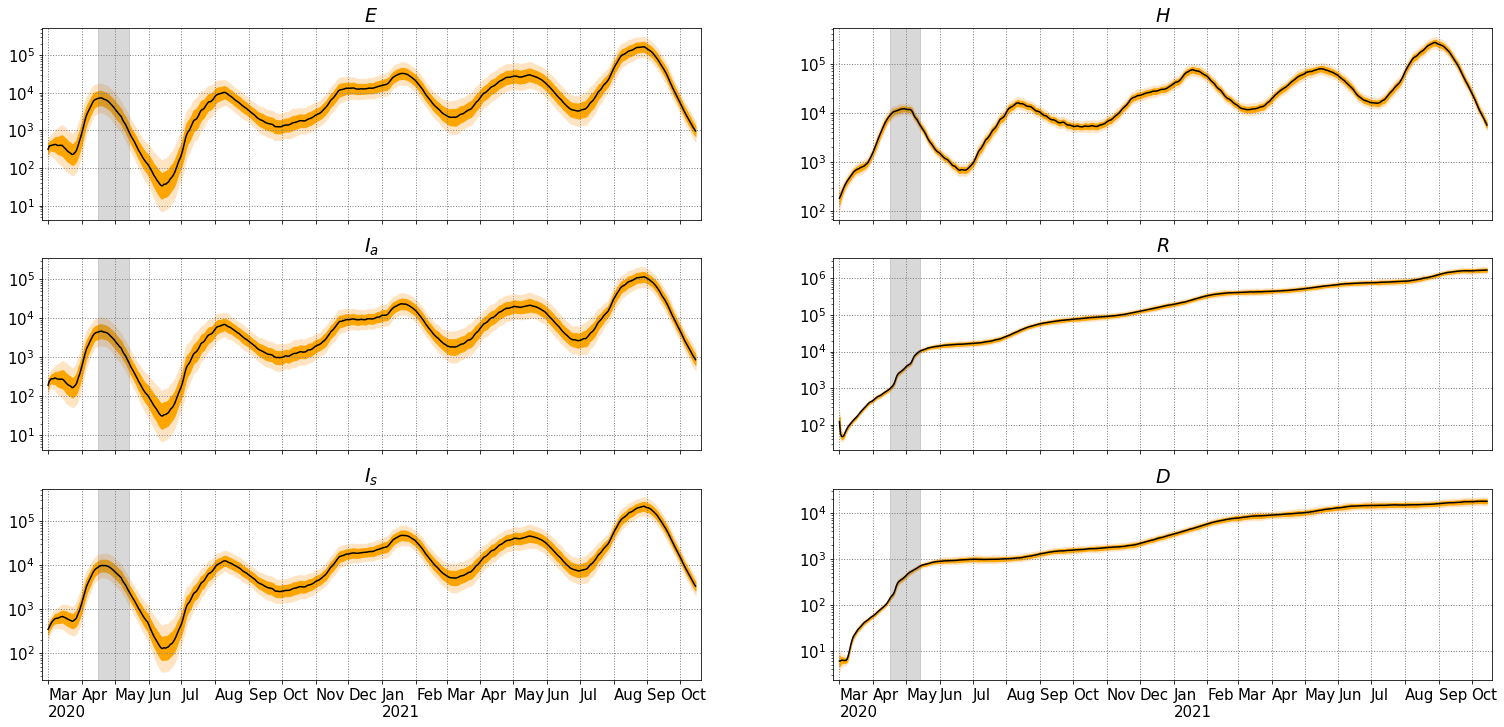}}
\caption{Analysis results of Japan}
\label{fig:com_ja}
\end{figure}

In Figure \ref{fig:br_o} we provide the analysis values of $\rt$ with $68\%$ CI and $95\%$ CI for each region. The additional red curves are again $\rt$ given by Toyokeizai.net. 
The  grey shaded regions mark the respective periods of state of emergency.
Though each region declared some states of emergency at different periods of time, 
they all have the common feature that the recent declarations are less efficient than 
the first ones. 
In Figure \ref{fig:com_os}, \ref{fig:com_ka}, and \ref{fig:com_ja}, we provide the analysis results for compartments $E$, $\is$, $\is$, $H$, $\rh$ and $\dh$ for each region. 

\begin{figure}[hbtp]
\centering
\subfigure{\includegraphics[width=15cm]{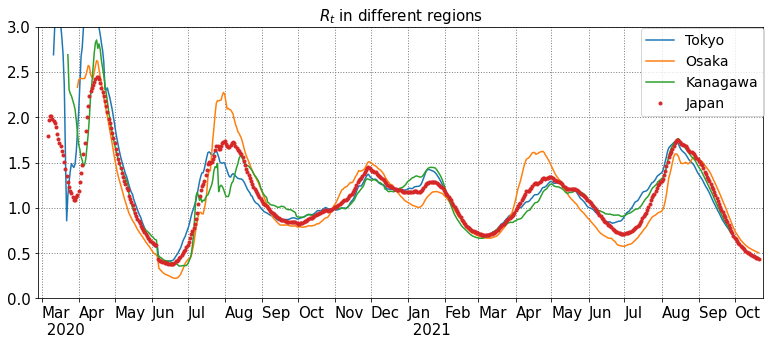}}
\caption{$\rt$ for each region}
\label{fig:rt_rg}
\end{figure}

In Figure \ref{fig:rt_rg}, the three analysis $\rt$ curves for Tokyo, Osaka and Kanagawa are drawn simultaneously.
The analysis $\rt$ for Japan is 
represented by red dots.
One observes that the general trends for each region are similar, but some delays are
also visible. For example in early April 2020, 
$\rt$ in Tokyo reached the first peak then started to decline, while $\rt$ in Kanagawa reached the first peak in the middle of April, and then start to decline. 
Subsequently, these two regions are quite correlated, but Osaka has a slightly independent behavior, both in the summer 2020 and in April 2021.

Based on these observations, one can suspect that the movement of populations might be related to the evolution of $\rt$ for regions which are close to each other. 
The study of the effect of movements can be carried out by constructing a proper model with multi-populations connected to each others. 
Some parameters may be tuned, for example for simulating movement restriction's policies. 
Such simulations could be useful for public health officials or for local governments to understand and estimate different disease-control strategies.
We plan to work in this direction in the future.

\section*{Acknowledgement}
S.~R. is supported by the grant \emph{Topological invariants through scattering theory and noncommutative geometry} from Nagoya University, and by JSPS Grant-in-Aid for scientific research C no 18K03328 \& 21K03292, and on leave of absence from Univ.~Lyon, Universit\'e Claude Bernard Lyon 1, CNRS UMR 5208, Institut Camille Jordan, 43 blvd. du 11 novembre 1918, F-69622 Villeurbanne cedex, France.


\begin{thebibliography}{10}

\bibitem{Arm}
E. Armstrong, M. Runge, J. Gerardin,
\emph{Identifying the measurements required to estimate rates of COVID-19 transmission, infection, and detection, using variational data assimilation},
Infectious Disease Modelling 6, 133--147, 2021.

\bibitem{AM}
F. Arroyo-Marioli, F. Bullano, S. Kucinskas,
C. Rond\' on-Moreno,
\emph{Tracking R of COVID-19: A new real-time
estimation using the Kalman filter},
PLoS ONE 16(1): e0244474, 2021.

\bibitem{BCB}
O. Byambasuren, M. Cardona, K. Bell, J. Clark, M.-L. McLaws, P. Glasziou, 
\emph{Estimating the extent of asymptomatic COVID-19 and its potential for community transmission: systematic review and meta-analysis},
J. Association of Medical Microbiology and Infectious Disease Canada
5 Issue 4, 223--234, 2020.

\bibitem{BDW}
F. Brauer, P.v.d. Driessche, J. Wu,
\emph{Mathematical epidemiology}, Lecture Notes in Mathematics 1945,
Springer, 2008.

\bibitem{BEC}
D. Buitrago-Garcia, D. Egli-Gany, M.J. Counotte, S. Hossmann, H. Imeri, A.M. Ipekci, et al.,  
\emph{Occurrence and transmission potential of asymptomatic and presymptomatic SARS-CoV-2 infections: A living systematic review and meta-analysis},  PLoS Med 17(9): e1003346, 2020.

\bibitem{BEM}
C.H. Bishop, B.J. Etherton, S.J. Majumdar, 
\emph{Adaptive sampling with the ensemble transform Kalman filter. Part I: Theoretical aspects}, 
Mon. Wea. Rev. 129, 420--436, 2001.

\bibitem{McA}
C. McAloon, A. Collins, K. Hunt, et al.,
\emph{Incubation period of COVID-19: a rapid systematic review and meta-analysis
of observational research}, BMJ Open 10:e039652, 2020.

\bibitem{DW}
P.v.d. Driessche, J. Watmough,
\emph{Reproduction numbers and sub-threshold endemic equilibria for compartmental models of disease transmission},
Mathematical Biosciences 180, 29--48, 2002.

\bibitem{E}
G. Evensen, J. Amezcua, M. Bocquet, A. Carrassi, A. Farchi, A. Fowler, P. L. Houtekamer, C. K. Jones, R. J. de Moraes, M. Pulido, C. Sampson, F. C. Vossepoel,
\emph{An international initiative of predicting the SARS-CoV-2 pandemic using ensemble data assimilation},  Foundations of Data Science,
American Institute of Mathematical Sciences, 2020.

\bibitem{E1}
G. Evensen, 
\emph{Sequential data assimilation with a nonlinear quasigeostrophic model using Monte Carlo methods to forecast error statistics},
J. Geophys. Res. 99, 10143--10162, 1994.

\bibitem{E2}
G. Evensen, 
\emph{The ensemble Kalman filter: Theoretical formulation and practical implementation}, 
Ocean Dynam. 53, 343--367, 2003.

\bibitem{E3}
G. Evensen, 
\emph{Data Assimilation: The Ensemble Kalman Filter}, 
Springer, 2006.

\bibitem{EN}
R. Engbert, M.M. Rabe, R. Kliegl, S. Reich,
\emph{Sequential Data Assimilation of the Stochastic SEIR
Epidemic Model for Regional COVID-19 Dynamics},
Bulletin of Mathematical Biology 83:1, 2021.

\bibitem{GH}
R. Ghostine, M. Gharamti, S. Hassrouny, I. Hoteit,
\emph{ An Extended SEIR Model with Vaccination for
Forecasting the COVID-19 Pandemic in Saudi Arabia Using an Ensemble
Kalman Filter},
Mathematics 9, 636. 2021.

\bibitem{GMc}
K.M. Gostic, L. McGough , E.B. Baskerville, S. Abbott, K. Joshi, C. Tedijanto, et al.,
\emph{Practical considerations for measuring the effective
reproductive number, $R_t$},
PLoS Comput Biol 16 No. 12,  e1008409, 21 pages, 2020.

\bibitem{H}
H.W. Hethcote,
\emph{The mathematics of infectious diseases}, 
SIAM Rev. 42 No. 4, 599--653, 2000.

\bibitem{HKS}
B.R. Hunt, E.J. Kostelich, I. Szunyogh, 
\emph{Efficient data assimilation for spatiotemporal chaos: a local ensemble transform Kalman filter}, 
Physica D. 230, 112--126, 2007.

\bibitem{HZ}
P.L. Houtekamer, F. Zhang, 
\emph{Review of the Ensemble Kalman filter for atmospheric data assimilation}, 
Mon. Wea. Rev. 144, 4489--4532, 2016.

\bibitem{K}
R.E. Kalman, 
\emph{A new approach to linear filtering and prediction problems},
Trans. ASME Ser. D: J. Basic Eng. 82, 35--45, 1960.

\bibitem{KB}
R.E. Kalman, R.S. Bucy, 
\emph{New results in linear filtering and prediction theory}, 
Trans. ASME Ser. D: J. Basic Eng. 83, 95--108, 1961.

\bibitem{KI}
T. Kuniya, H. Inaba,
\emph{Possible effects of mixed prevention strategy for COVID-19 epidemic: massive
testing, quarantine and social distancing}, 
AIMS Public Health 7 No 3, 490--503, 2020.

\bibitem{L}
M. Li, 
\emph{An Introduction to Mathematical Modeling of Infectious Diseases}, 
Mathematics of Planet Earth 2, Springer 2018.

\bibitem{LI}
R. Li, et al.,
\emph{Substantial undocumented infection
facilitates the rapid dissemination of novel
coronavirus (SARS-CoV-2)},
Science 368, 489–493, 2020.

\bibitem{M3}
m3: \url{https://www.m3.com/open/iryoIshin/article/849820/}

\bibitem{MI}
L. Mitchell, A. Arnold,
\emph{Analyzing the effects of observation function selection in ensemble Kalman filtering for epidemic models},
Mathematical Biosciences 339, 2021.

\bibitem{MLIT}
MLIT: \url{https://www.mlit.go.jp/tetudo/tetudo\_fr1\_000062.html}

\bibitem{NA}
P. Nadler,  S. Wang,  R. Arcucci, X. Yang, Y. Guo, 
\emph{An epidemiological modelling approach for COVID-19 via data
assimilation},
European Journal of Epidemiology 35, 749--761, 2020.

\bibitem{NH}
H. Nishiura: \url{https://github.com/contactmodel/COVID19-Japan-Reff}

\bibitem{OHS1}
E. Ott, B.R. Hunt, I. Szunyogh, M. Corazza, E. Kalnay, D.J. Patil, J.A. Yorke, A.V. Zimin, E.J. Kostelich, 
\emph{Exploiting local low dimensionality of the atmospheric dynamics for efficient ensemble Kalman filtering},
Preprint: \url{https://arxiv.org/abs/physics/0203058v3}.

\bibitem{OHS2}
E. Ott, B.R. Hunt, I. Szunyogh, A.V. Zimin, E.J. Kostelich, M. Corazza, E. Kalnay, D.J. Patil, J.A. Yorke, 
\emph{A local ensemble Kalman filter for atmospheric data assimilation}, 
Tellus A 56, 415--428, 2004.

\bibitem{OP}
Osaka prefecture government, 
\emph{Citizens awareness and behavior change of measures against COVID-19}, 
\url{http://www.pref.osaka.lg.jp/hodo/attach/hodo-40479\_4.pdf}

\bibitem{PL}
A.M. Pollock, J. Lancaster, 
\emph{Asymptomatic transmission of covid-19}, 
BMJ 371:m4851, 2020. 

\bibitem{RE}
T.C. Rebollo,  D. Coronil,
\emph{Predictive data assimilation through Reduced Order Modeling for epidemics with data uncertainty},
Preprint: \url{https://arxiv.org/abs/2004.12341}.

\bibitem{RH}
C. J. Rhodes, T. D. Hollingsworth,
\emph{Variational data assimilation with epidemic models},
Journal of Theoretical Biology 258, 591--602, 2009.

\bibitem{SRM}
C. Sun, S. Richard, T. Miyoshi,
\emph{Agent-based model and data assimilation: 
Analysis of COVID-19 in Tokyo},
Preprint: \url{https://arxiv.org/abs/2109.00258}.

\bibitem{TAB}
M.K. Tippett, J.L. Anderson, C.H. Bishop, T.M. Hamill, J.S. Whitaker,
\emph{Ensemble square-root filters}, 
Mon. Wea. Rev. 131, 1485--1490, 2003.

\bibitem{TBP}
Bureau of social welfare and public health, 
\emph{About death cases due to COVID-19 in Tokyo},
\url{https://www.fukushihoken.metro.tokyo.lg.jp}

\bibitem{TC}
Tokyo metropolitan government, 
\emph{COVID-19 The information website},
\url{https://stopcovid19.metro.tokyo.lg.jp}.

\bibitem{TK}
Toyokeizai,
\emph{Coronavirus disease (COVID-19) situation report in Japan},
\url{https://toyokeizai.net/sp/visual/tko/covid19/en.html}.

\bibitem{WH}
J.S. Whitaker, T.M. Hamill, 
\emph{Ensemble data assimilation without perturbed observations}, 
Mon. Wea. Rev. 130, 1913--1924, 2002.

\bibitem{WHO}
World Health Organization:
\url{https://www.who.int/}.

\end{thebibliography}
\end{document}